\documentclass[10pt]{article}
\usepackage{mypackages}
\usepackage{mycommands}

\title{Explicit stabilized multirate method for stiff stochastic differential equations}

\author{Assyr Abdulle\thanks{\email{assyr.abdulle@epfl.ch}} \and Giacomo Rosilho de Souza\thanks{\email{giacomo.rosilhodesouza@epfl.ch}}}

\begin{document}

\maketitle

\blfootnote{\hspace{-4ex}\'Ecole Polytechnique F\'ed\'erale de Lausanne (EPFL), SB-MATH-ANMC, Station 8, 1015 Lausanne, Switzerland.}

\begin{abstract}
Stabilized explicit methods are particularly efficient for large systems of stiff stochastic differential equations (SDEs) due to their extended stability domain. However, they lose their efficiency when a severe stiffness is induced by very few ``fast'' degrees of freedom, as the stiff and nonstiff terms are evaluated concurrently. Therefore, inspired by 
[A.~Abdulle, M.~J. Grote, and G.~{Rosilho de Souza}, \newblock {\em Preprint} (2020), \burlalt{arXiv:2006.00744}{http://arxiv.org/abs/2006.00744}],
we introduce a stochastic modified equation whose stiffness depends solely on the ``slow'' terms.
By integrating this modified equation with a stabilized explicit scheme we devise a multirate method which overcomes the bottleneck caused by a few severely stiff terms and recovers the efficiency of stabilized schemes for large systems of nonlinear SDEs. The scheme is not based on any scale separation assumption of the SDE and therefore it is employable for problems stemming from the spatial discretization of stochastic parabolic partial differential equations on locally refined grids. The multirate scheme has strong order $1/2$, weak order $1$ and its stability is proved on a model problem. Numerical experiments confirm the efficiency and accuracy of the scheme.
\end{abstract}


\section{Introduction}\label{sec:intro}
We consider Itô systems of stochastic differential equations of the form
%\begin{subequations}%\label{eq:msde}
\begin{equation}\label{eq:msde}
\dif X(t) = f(X(t))\dif t + g(X(t))\dif W(t), \qquad\qquad X(0)=X_0,
\end{equation}
where
\begin{equation}%\label{eq:msdeb}
f(X)=\ff(X)+\fs(X)
\end{equation}
%\end{subequations}
splits in an inexpensive but stiff term $\ff$ associated to fast time-scales and an expensive but mildly stiff term $\fs$ associated to relatively slow time-scales. In \eqref{eq:msde}, $X(t)$ is a stochastic process in $\Rn$, $\ff,\fs:\Rn\rightarrow\Rn$ are drift terms, $g:\Rn\rightarrow\Rb^{n\times \modr{l}}$ is the diffusion term and $W(t)$ is an $\modr{l}$-dimensional Wiener process. 
We emphasize that $\ff$ is stiff compared to $\fs$, nonetheless not all the eigenvalues of the Jacobian of $\ff$ are large in magnitude, hence \emph{we do not make any scale separation assumption}. Therefore, the schemes presented here can be employed, for instance, for problems stemming from the spatial discretization of stochastic parabolic partial differential equations on locally refined grids. Indeed, $\ff$ and $\fs$ would represent the discrete Laplacian in the refined and coarse region, respectively; hence, $\ff$ contains \emph{fast and slow} scales. In contrast, $\fs$ contains relatively slow terms only.

Due to the stiffness of $f$, traditional explicit schemes as Euler--Maruyama face stringent conditions on the step size. On the other hand, implicit methods require the solution to possibly nonlinear systems. Stochastic stabilized explicit methods (the S-ROCK family) \cite{AAV18,AbC07,AbC08,AVZ13b} are a good compromise, as they enjoy an extended stability domain growing \textit{quadratically} with the number of stages $s$. The scheme presented in \cite{AAV18}, called SK-ROCK for second kind Runge--Kutta orthogonal Chebyshev, attains an optimal mean-square stability domain of size $L_s\approx 2 s^2$. It is based on the deterministic Runge--Kutta--Chebyshev (RKC) method \cite{SSV98,HoS80,VHS90}, which has an optimal stability domain along the negative real axis for an $s$-stage Runge--Kutta method \cite{GuL60}, and employs second kind Chebyshev polynomials for the stabilization of the stochastic integral. However, the number of stages $s$ is dictated by the stiffness of $f$. Therefore, even if stiffness is induced by only a few degrees of freedom in $\ff$, the cost of numerical integration is high; indeed, the nonstiff expensive term $\fs$ is evaluated concurrently to the stiff term $\ff$. Consequently, a multirate/multiscale strategy must be employed.

In the class of multiscale methods for stochastic differential equations, we find the heterogeneous multiscale methods \cite{ELV05a,GKK06,Liu10,Van03}. They are based on a scale separation assumption and therefore derive an effective equation for the slow variables, which depends on the invariant measure of the fast dynamics. 
An extension of those methods to stochastic partial differential equations is found in \cite{AbP12,APV17}, while a close family of schemes are the projective methods \cite{GKK06,PaK07} --- see \cite{Van07} for a review. As the aforementioned methods are strongly based on a scale separation assumption, they cannot be employed when \eqref{eq:msde} stems from the spatial discretization of a stochastic parabolic partial differential equation.

Since the early work of Rice \cite{Ric60}, many multirate strategies for the solution of the stiff ordinary differential equation (ODE) $y'=\ff(y)+\fs(y)$ have been developed, see for instance \cite{And79,EnL97,GeW84,GKR01,Kva99,SHV07,SkA89}.
%\cite{And79,EnL97,GeW84,GKR01,GuR93,Hof76,Kva99,SHV07,SaM10,SkA89}
These methods are based on predictor-corrector strategies, on interpolation/extrapolation of ``fast'' and ``slow'' variables (which is known to trigger instabilities) or are implicit. An alternative approach consists in deriving an effective equation for the slow dynamics \cite{E03,EnT05,GIK03}, but this strategy works for scale separated problems only. More recently, multirate methods based on the GARK framework have been developed \cite{GuS16,RSS20,San19,SaG15}. This approach allows for the development of high order multirate schemes but in order to obtain satisfying stability properties some degree of implicitness is required.
In \cite{AGR20}, a stabilized explicit multirate method, called $\mRKC$ for multirate RKC, is introduced. It is based on a \emph{modified equation}, defined by an \emph{averaged force}, whose stiffness depends on $\fs$ only and is decreased due to an average along the direction defined by a fast but cheap \emph{auxiliary problem}. Due to the decreased stiffness, integration of the modified equation by an explicit scheme is cheaper than integrating the original problem with the same scheme. In \cite{AGR20}, the modified equation and the auxiliary problems are integrated by RKC schemes; the number of expensive evaluations of $\fs$ depends on the slow terms only and the bottleneck caused by the stiffness of $\ff$ is overcome without sacrificing accuracy nor explicitness. 

The contribution of this paper is twofold. First, in \cref{sec:modeq} we extend the modified equation for ODEs, introduced in \cite{AGR20}, to SDEs, obtaining a \emph{stochastic modified equation}. This is not a trivial generalization of \cite{AGR20} as it requires an approximation of the diffusion term $g$, called \emph{damped diffusion}, so that the mean-square stability properties of \cref{eq:msde} are inherited by the stochastic modified equation. Second, in \cref{sec:mskrock} we define the multirate SK-ROCK ($\mSKROCK$) method as a time discretization of the stochastic modified equation using the SK-ROCK scheme, while the deterministic auxiliary problems are solved with RKC schemes. The resulting method inherits the main properties of the $\mRKC$ and SK-ROCK schemes: it is explicit, the stability domain grows optimally and quadratically with the number of stages, the number of expensive function evaluations depends on $\fs$ only, it is not based on any scale separation assumption, there is no need of interpolations nor extrapolations and therefore it is straightforward to implement. The stability and accuracy analysis of the mSK-ROCK scheme is presented in \cref{sec:analysis}, while \cref{sec:numexp} is devoted to numerical experiments, where we illustrate the theoretical results and confirm the efficiency of the multirate method. Application of the method to the E. Coli bacteria heat shock response and to a diffusion problem across a narrow channel with multiplicative time-space noise is also provided.

%\section{The multirate SK-ROCK method}\label{sec:mskrock}
\section{The stochastic modified equation}\label{sec:modeq}
In this section we introduce the \emph{stochastic modified equation}
\begin{equation}\label{eq:modsde}
\dif \Xe(t) = \fe(X(t))\dif t +\ge(X(t))\dif W(t), \qquad\qquad  \Xe(0)=X_0,
\end{equation}
which is an approximation of \cref{eq:msde} but whose stiffness depends solely on $\fs$ and therefore is not affected by the severely stiff terms in $\ff$. 
Indeed, the \emph{averaged force} $\fe$ is an approximation of $f$ satisfying $\rhoe\leq\rhos$, where $\rhoe,\rhos$ are the spectral radii of the Jacobians of $\fe,\fs$, respectively. Hence, in \cref{eq:modsde}, the aim in replacing $\ff+\fs$ by $\fe$ is to reduce the stiffness. Differently, the \emph{damped diffusion} $\ge$ is needed to preserve the mean-square stability properties of the original problem \cref{eq:msde}; as $\fe$ is less stiff than $f$ it is also less contractive and therefore it cannot damp the original noise term $g$ enough to maintain stability. \modr{Note that mean-square stability implies that the expectation of the squared norm of perturbations vanish at infinity, for \cref{eq:modsde} it is studied on a model problem in \cref{sec:modeqstoc} below.}

We first recall the averaged force and the deterministic modified equation introduced in \cite{AGR20}. Then we define the damped diffusion $\ge$ and analyze the stochastic modified equation \cref{eq:modsde}. \modr{In order to ensure existence and uniqueness of the solutions to \cref{eq:msde} and the next \cref{eq:mrode,eq:defu,eq:defv}, here and in the foregoing sections we assume that $\ff,\fs$ and $g$ are uniformly Lipschitz continuous and satisfy a linear growth condition. Under the same assumptions, with a few tedious computations it is possible to show that $\fe,\ge$ defined below are also uniformly Lipschitz continuous with linear growth and therefore the solutions to \cref{eq:modsde,eq:modode} also exist and are unique. See \cite[Lemma 4.16]{Ros20}.}

\subsection{The modified equation for deterministic problems}
We consider stiff multirate differential equations of the type
\begin{equation}\label{eq:mrode}
y' =f(y):= \ff(y)+\fs(y), \qquad\qquad y(0)=y_0,
\end{equation}
where $\ff$ is a cheap but severely stiff term and $\fs$ is an expensive but only mildly stiff term. In \cite{AGR20}, the right-hand side $f=\ff+\fs$ is replaced by an averaged force $\fe$ depending on a free parameter $\eta\geq 0$. For large enough $\eta$ it holds $\rhoe\leq \rhos$ and since $\rhos\ll \rhof$, where $\rhof,\rhos$ are the spectral radii of the Jacobians of $\ff,\fs$, respectively, integration of the averaged system
\begin{equation}\label{eq:modode}
\ye'=\fe(\ye), \qquad\qquad  \ye(0)=y_0
\end{equation}
with an explicit scheme is much cheaper than \cref {eq:mrode}. In practice, evaluation of $\fe$ requires the solution to a fast but cheap auxiliary ODE, which is as well approximated by an explicit scheme. In the rest of the section we will define \eqref{eq:modode} and recall some its key stability properties.

\subsubsection*{The averaged force}
\begin{figure}
	\centering
	\begin{tikzpicture}[scale=\plotimscale]
	\begin{axis}[height=\aspectratio*\plotimsizeu\textwidth,width=\plotimsizeu\textwidth, ymin=-0.2, ymax=1,xmax=0,xmin=-100,legend columns=1,legend style={draw=\legendboxdraw,fill=\legendboxfill,at={(0,1)},anchor=north west},
	xlabel={$z$}, ylabel={},label style={font=\normalsize},tick label style={font=\normalsize},legend image post style={scale=\legendmarkscale},legend style={nodes={scale=\legendfontscale, transform shape}},grid=none]
	\addplot[color=colorone,line width=\plotlinewidth pt,mark=\markone,mark repeat=20,mark phase=0,mark size=\plotmarksizeu pt] table [x=z,y=phiz,col sep=comma] 
	{data/text/phis.csv};\addlegendentry{$\varphi(z)$}
	\addplot[color=colortwo,line width=\plotlinewidth pt,mark=\marktwo,mark repeat=20,mark phase=10,mark size=\plotmarksizeu pt] table [x=z,y=phizd2sq,col sep=comma] 
	{data/text/phis.csv};\addlegendentry{$\varphi(z/2)^2$}
	\end{axis}
	\end{tikzpicture}
	\caption{Illustration of $\varphi(z)$ and $\varphi(z/2)^2$, used for damping the drift and diffusion terms, respectively.}
	\label{fig:phi}
\end{figure}
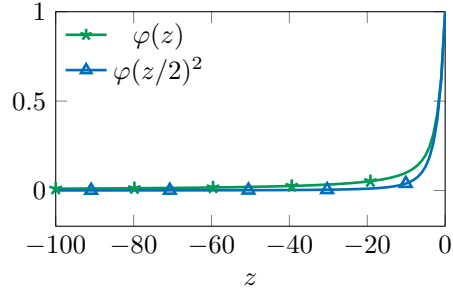

Here we define the averaged force $\fe$ and recall some of its main properties.
\begin{definition}\label{def:fe}
	\modr{For $\eta>0$, the averaged force $\fe:\Rn\rightarrow\Rn$ is defined as}
	\begin{equation}\label{eq:deffedif}
		\fe(y)=\frac{1}{\eta}(u(\eta)-y),
	\end{equation}
	\modr{where {\it the auxiliary solution} $u:[0,\eta]\rightarrow \Rn$ is defined by {\it the auxiliary equation}}
	\begin{align}\label{eq:defu}
		u'&=\ff(u)+\fs(y), & 
		u(0)=y.
	\end{align}
	\modr{For $\eta=0$, let $f_0=f$ (note that $f_0=\lim_{\eta\to 0^+}\fe$).}
\end{definition}
\modr{Hence, in \cref{eq:modode} whenever $\fe(\ye(t))$ is evaluated the auxiliary problem \cref{eq:defu} is solved with initial value $u(0)=\ye(t)$.}
From \cref{eq:defu,eq:deffedif} we obtain
\begin{equation}\label{eq:deffeint}
\fe(y)=\frac{1}{\eta}\int_0^\eta u'(s)\dif s =\fs(y)+\frac{1}{\eta}\int_0^\eta \ff(u(s))\dif s,
\end{equation}
hence $\fe$ is an average of $f$ along the \emph{auxiliary solution} $u$. In \cite{AGR20} it is shown that $\ff$ has a smoothing effect on $\fe$ and thus \eqref{eq:modode} has a reduced stiffness when compared to \eqref{eq:mrode}. More precisely, for linear $\ff$ we have the following result.
\begin{lemma}\label{lemma:phif}
	Let $\ff(y)=A_F\,y$ with $A_F\in\Rb^{n\times n}$. Then
	\begin{equation}\label{eq:deffephi}
	\fe(y)=\varphi(\eta A_F)f(y),
	\end{equation}
where 
\begin{equation}\label{eq:defphi}
	\varphi(z)=\frac{e^z-1}{z} \;\mbox{ for }\;z\neq 0\qquad \mbox{and}\qquad \varphi(0)=1.
\end{equation}
%and $\varphi(0)=1$ is defined by continuous extension.
\end{lemma}
The function $\varphi(z)$ satisfies $\lim_{z\to -\infty}\varphi(z)=0$ and $\varphi(z)\in (0,1)$ for all $z<0$, see \cref{fig:phi}. \modr{Note that $\varphi(z)$ is an entire function and therefore its evaluation on square matrices is well defined.} \Cref{lemma:phif} states that if $\ff$ is linear then we have a closed expression for $\fe$ and we see in \cref{eq:deffephi} the smoothing effect of a negative definite matrix $A_F$ on $f$. In \cref{eq:deffephi} we see as well the role of $\eta$: it is a free parameter used to tune this smoothing effect. In \cite{AGR20} it is shown that $\fe$ inherits the contractivity properties of $f$ and that the error between the exact solution $y$ to \cref{eq:mrode} and the solution $\ye$ to \cref{eq:modode} is of first-order in $\eta$ and bounded independently of the stiffness of the problem.

\subsubsection*{Linear stability analysis on the multirate test equation}\label{sec:stab_an_modode}
Here we recall the conditions on $\eta$ for which the spectral radius of $\fe$ depends only on $\fs$. To do so, we apply \cref{def:fe} to the \emph{multirate test equation}
\begin{equation}\label{eq:mrtesteq}
y' = \lambda y+\zeta y, \qquad\qquad y(0)=y_0,
\end{equation}
with $\lambda,\zeta\leq 0$ and $y_0\in\Rb$. We set $\ff(y)=\lambda y$ and $\fs(y)=\zeta y$; thus, $\rhof=|\lambda|$ and $\rhos=|\zeta|$. From \cref{eq:deffephi} follows
\begin{equation}\label{eq:phiefetest}
\fe(y)=\varphi(\eta\lambda)(\lambda+\zeta)y
\end{equation}
and thus \cref{eq:modode} becomes
\begin{equation}\label{eq:linscamod}
\ye'=\varphi(\eta\lambda)(\lambda+\zeta)\ye, \qquad\qquad \ye(0)=y_0.
\end{equation}
The next \cref{thm:eta}, proved in \modr{\cite[Theorem 2.7]{AGR20}}, states that if $\eta$ is taken large enough then $|\varphi(\eta\lambda)(\lambda+\zeta)|\leq |\zeta|$ and the stiffness of \cref{eq:linscamod} depends only on $\zeta$, thus on $\fs$. Furthermore, $\lambda$ can take any nonpositive value and thus there is no scale separation assumption.
\begin{theorem}\label{thm:eta}
Let $\zeta<0$, it holds $\varphi(\eta\lambda)(\lambda+\zeta)\in [\zeta,0]$ for all $\lambda\leq 0$ if, and only if, $\eta|\zeta|\geq 2$.
\end{theorem}

\subsection{The modified equation for stochastic problems}\label{sec:modeqstoc}
As $\fe$ is less stiff than $f$ it follows that it has also weaker contractivity and therefore it cannot control the original noise term $g$, hence in this section we introduce a damped noise term $\ge$ to restore for the modified equation \cref{eq:modsde} the mean-square stability properties of the original problem \cref{eq:msde}.

\subsubsection*{The damped diffusion}\label{sec:defge}
Here we define the damped diffusion term $\ge$ of \eqref{eq:modsde} and study its properties. \modr{We consider here a vector valued diffusion term $g:\Rn\rightarrow\Rn$ for simplicity. For a matrix valued diffusion $g:\Rn\rightarrow\Rb^{n\times l}$ we can simply apply the same definitions and results column-wise. In \cref{sec:mskrockalg} we will also discuss how to preserve the scheme efficiency disregarding the number of columns in $g$.}
\begin{definition}\label{def:ge}
\modr{Let $\eta> 0$, the damped diffusion $\ge:\Rn\rightarrow\Rn$ is defined as}
\begin{equation}\label{eq:defgedif}
	\ge(x)=\frac{1}{\eta}(v(\eta)-\bv (\eta)),
\end{equation}
\modr{where the auxiliary solutions $v,\bv :[0,\eta]\rightarrow\Rn$ are defined by the auxiliary equations}
\begin{equation}\label{eq:defv}
v'=\frac{1}{2}\ff(v)+g(x), \qquad\qquad \bv '=\frac{1}{2}\ff(\bv ), \qquad\qquad v(0)=\bv(0)=x.
\end{equation}
\modr{For $\eta=0$ let $g_0=g$.}
\end{definition}
The motivation for \cref{def:ge} and the factor $1/2$ will be better seen in the linear stability analysis given in the next paragraph.

\begin{lemma}\label{lemma:phig}
	Let $\ff(y)=A_F\,y$ with $A_F\in\Rb^{n\times n}$, then
	\begin{equation}\label{eq:defgephi}
	\ge(x)=\varphi\left(\frac{\eta}{2}A_F\right)g(x).
	\end{equation}
\end{lemma}
\begin{proof}
	Replacing $\ff(y)=A_F\, y$ in \eqref{eq:defv} and using the variation-of-constants formula we deduce
	\begin{equation}
	v(\eta)=e^{\frac{\eta}{2}A_F}x+\eta\varphi\left(\frac{\eta}{2}A_F\right)g(x), \qquad\qquad \bv (\eta)=e^{\frac{\eta}{2}A_F}x
	\end{equation}
	and \eqref{eq:defgephi} follows from \eqref{eq:defgedif}.
\end{proof}
In \eqref{eq:defgephi} we observe the smoothing effect of $\ff$ on $g$ and since $\varphi(z)=1+\bigo{z}$ as $z\to 0$ then $\ge(x)=g(x)+\bigo{\eta}$ as $\eta\to 0$.
For a general $\ff$, from \cref{eq:defv}, \cref{eq:defgedif}, we obtain
\begin{equation}\label{eq:defgeint}
	\ge(x)=\frac{1}{\eta}\int_0^\eta (v'(s)-\bv '(s))\dif s = g(x)+\frac{1}{2\eta}\int_0^\eta ( \ff(v(s))-\ff(\bv (s)))\dif s,
\end{equation}
hence $\ge$ is still composed of $g$ plus additional higher order terms. The role of $\ff(v)$ is still to stabilize $g$, while $\ff(\bv )$ is used to remove the low order polluting terms introduced by $\ff(v)$ (as is seen in the proof of \cref{lemma:phig}).

\subsubsection*{Linear mean-square stability analysis of the modified equation}\label{sec.stab_an_modsde}
As for stiff SDEs, we will consider the relevant notion of mean-square stability and extend the widely used linear scalar test equation \cite{Hig01,SaM96} to multirate stochastic problems. Therefore, we consider
\begin{equation}\label{eq:smrtesteq}
\dif X(t) = (\lambda +\zeta) X(t)\dif t +\mu X(t)\dif W(t), \qquad\qquad X(0)=X_0,
\end{equation}
with $\lambda,\zeta\leq 0$ and $\mu\in\Rb$. Next, we identify $\ff(X)=\lambda X$ and $\fs(X)=\zeta X$, while we let $g(X)=\mu X$.
The exact solution to \cref{eq:smrtesteq} is called mean-square stable if, and only if, $\lim_{t\to\infty} \exp(|X(t)|^2)=0$, which holds if $(\lambda,\zeta,\mu)\in \mathcal{S}^{mMS}$, where
\begin{equation}\label{eq:defSmMS}
\mathcal{S}^{mMS} = \{ (\lambda,\zeta,\mu)\in \Rb^3 \, :\, \lambda+\zeta+\frac{1}{2}|\mu|^2<0, \lambda\leq 0, \zeta\leq 0 \}
\end{equation}
is the mean-square stability domain for the stochastic multirate test equation \cref{eq:smrtesteq}. 

From \cref{eq:deffephi,eq:defgephi} the modified equation \cref{eq:modsde} yields for the stochastic multirate test equation 
\begin{equation}\label{eq:smrtesteqmod}
\dif \Xe(t) = \varphi(\eta\lambda)(\lambda+\zeta) \Xe(t)\dif t +\varphi\left(\frac{\eta}{2}\lambda\right)\mu \Xe(t)\dif W(t),\qquad\qquad \Xe(0)=X_0.
\end{equation}
In \cref{thm:msstab_modsde} we will show that \cref{eq:smrtesteqmod} is mean-square stable, to do so the next property of $\varphi(z)$ (defined in \cref{eq:defphi}) is crucial. 
\begin{lemma}\label{lemma:phisq}
Let $z\in\Rb$, then $\varphi\left(\frac{z}{2}\right)^2\leq \varphi(z)$.
\end{lemma}
\begin{proof}
	Since $\varphi(z)=\int_0^1 e^{z s}\dif s$ the result follows from Jensen's inequality. Indeed,
	\begin{equation}
	\varphi\left(\frac{z}{2}\right)^2=\left(\int_0^1 e^{\frac{z}{2}s}\dif s\right)^2\leq\int_0^1 e^{z s}\dif s =\varphi(z).
	\tag*{\qedhere}
	\end{equation}
\end{proof}

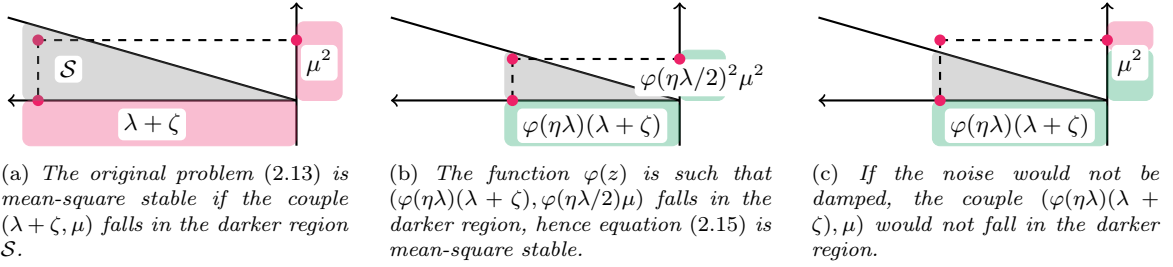
\begin{figure}
	\begin{center}
		\begin{subfigure}[t]{0.3\textwidth}
			\centering
			\begin{tikzpicture}[scale=\plotimscale]
			\coordinate (a) at (0,0);
			\coordinate (b) at (3.8,0);
			\coordinate (c) at (3.8,1.3);
			\coordinate (d) at (3.8,-0.6);
			\coordinate (e) at (1.9,-0.1);
			\coordinate (f) at (0,1.1);
			\coordinate (g) at (0.8,0.6);
			\coordinate (r1) at (0.2,-0.6);
			\coordinate (r2) at (0.2,0);
			\coordinate (r3) at (0.2,1.042);
			\coordinate (t1) at (4.4,1.042);
			\coordinate (t2) at (4.1,0.5);
			\draw[thick,<-] (a)--(b);
			\draw[thick,<-] (c)--(d);
			\draw[thick,-] (b)--(f);
			\draw[thick,dashed,-] (0.4,0.8)--(3.8,0.8) node[circle,fill,inner sep=1.5pt,colorthree]{}; 
			\draw[thick,dashed,-] (0.4,0.8) node[circle,fill,inner sep=1.5pt,colorthree]{}--(0.4,0) node[circle,fill,inner sep=1.5pt,colorthree]{}; 
			\draw[colorthree,fill=colorthree,rounded corners=1mm,opacity=0.3] (r1) rectangle (b);
			\draw[colorthree,fill=colorthree,rounded corners=1mm,opacity=0.3] (t1) rectangle (b);
			\draw[gray,fill=gray,rounded corners=1mm,opacity=0.3] (r2)--(r3)--(b)--cycle;
			\node[fill=white,rounded corners=2pt,inner sep=2pt,anchor=north] at (g) {{\small $\mathcal{S}$}};
			\node[fill=white,rounded corners=2pt,inner sep=2pt,anchor=north] at (e) {{\small $\lambda+\zeta$}};
			\node[fill=white,rounded corners=2pt,inner sep=2pt,anchor=center] at (t2) {{\small $\mu^2$}};
			\end{tikzpicture}
			\caption{The original problem \cref{eq:smrtesteq} is mean-square stable if the couple $(\lambda+\zeta,\mu)$ falls in the darker region $\mathcal S$.}
			\label{fig:origmsde}
		\end{subfigure}\hfill 
		\begin{subfigure}[t]{0.33\textwidth}
			\centering
			\begin{tikzpicture}[scale=\plotimscale]
			\coordinate (a) at (0,0);
			\coordinate (b) at (3.8,0);
			\coordinate (c) at (3.8,1.3);
			\coordinate (d) at (3.8,-0.6);
			\coordinate (e) at (2.65,-0.1);
			\coordinate (f) at (0,1.1);
			\coordinate (g) at (0.8,0.6);
			\coordinate (r1) at (1.5,-0.6);
			\coordinate (r2) at (1.5,0);
			\coordinate (r3) at (1.5,0.666);
			\coordinate (t1) at (4.4,0.666);
			\coordinate (t2) at (4.1,0.3);
			\draw[thick,<-] (a)--(b);
			\draw[thick,<-] (c)--(d);
			\draw[thick,-] (b)--(f);
			\draw[colorone,fill=colorone,rounded corners=1mm,opacity=0.3] (r1) rectangle (b);
			\draw[colorone,fill=colorone,rounded corners=1mm,opacity=0.3] (t1) rectangle (b);
			\draw[gray,fill=gray,rounded corners=1mm,opacity=0.3] (r2)--(r3)--(b)--cycle;
%			\node[fill=white,rounded corners=2pt,inner sep=2pt,anchor=north] at (g) {{\small $\mathcal{S}$}};
			\node[fill=white,rounded corners=2pt,inner sep=2pt,anchor=north] at (e) {{\small $\varphi(\eta\lambda)(\lambda+\zeta)$}};
			\node[fill=white,rounded corners=2pt,inner sep=2pt,anchor=center] at (t2) {{\small $\varphi(\eta\lambda/2)^2\mu^2$}};
			\draw[thick,dashed,-] (1.6,0.55)--(3.8,0.55) node[circle,fill,inner sep=1.5pt,colorthree]{}; 
			\draw[thick,dashed,-] (1.6,0.55) node[circle,fill,inner sep=1.5pt,colorthree]{}--(1.6,0) node[circle,fill,inner sep=1.5pt,colorthree]{}; 
			\end{tikzpicture}
			\caption{The function $\varphi(z)$ is such that $(\varphi(\eta\lambda)(\lambda+\zeta),\varphi(\eta\lambda/2)\mu)$ falls in the darker region, hence equation \cref{eq:smrtesteqmod} is mean-square stable.}
			\label{fig:modsdea}
		\end{subfigure}\hfill 
		\begin{subfigure}[t]{0.3\textwidth}
			\centering
			\begin{tikzpicture}[scale=\plotimscale]
			\coordinate (a) at (0,0);
			\coordinate (b) at (3.8,0);
			\coordinate (c) at (3.8,1.3);
			\coordinate (d) at (3.8,-0.6);
			\coordinate (e) at (2.65,-0.1);
			\coordinate (f) at (0,1.1);
			\coordinate (g) at (0.8,0.6);
			\coordinate (r1) at (1.5,-0.6);
			\coordinate (r2) at (1.5,0);
			\coordinate (r3) at (1.5,0.666);
			\coordinate (t1) at (4.4,0.666);
			\coordinate (t2) at (4.1,0.5);
			\coordinate (t3) at (3.8,1.042);
			\draw[thick,<-] (a)--(b);
			\draw[thick,<-] (c)--(d);
			\draw[thick,-] (b)--(f);
			\draw[colorone,fill=colorone,rounded corners=1mm,opacity=0.3] (r1) rectangle (b);
			\draw[colorone,fill=colorone,rounded corners=1mm,opacity=0.3] (t1) rectangle (b);
			\draw[colorthree,fill=colorthree,rounded corners=1mm,opacity=0.3] (t1) rectangle (t3);
			\draw[gray,fill=gray,rounded corners=1mm,opacity=0.3] (r2)--(r3)--(b)--cycle;
			%			\node[fill=white,rounded corners=2pt,inner sep=2pt,anchor=north] at (g) {{\small $\mathcal{S}$}};
			\node[fill=white,rounded corners=2pt,inner sep=2pt,anchor=north] at (e) {{\small $\varphi(\eta\lambda)(\lambda+\zeta)$}};
			\node[fill=white,rounded corners=2pt,inner sep=2pt,anchor=center] at (t2) {{\small $\mu^2$}};
			\draw[thick,dashed,-] (1.6,0.8)--(3.8,0.8) node[circle,fill,inner sep=1.5pt,colorthree]{}; 
			\draw[thick,dashed,-] (1.6,0.8) node[circle,fill,inner sep=1.5pt,colorthree]{}--(1.6,0) node[circle,fill,inner sep=1.5pt,colorthree]{}; 
			\end{tikzpicture}
			\caption{If the noise would not be damped, the couple $(\varphi(\eta\lambda)(\lambda+\zeta),\mu)$ would not fall in the darker region.}
			\label{fig:modsdeb}
		\end{subfigure}
	\end{center}
	\caption{Illustration of the mean-square stability properties of the original problem \cref{eq:smrtesteq} (left), the stochastic modified equation \cref{eq:smrtesteqmod} (center) and a case where the noise is not damped (right).}
	\label{fig:stabmodsde}
\end{figure}

\begin{theorem}\label{thm:msstab_modsde}
If the stochastic multirate test equation \cref{eq:smrtesteq} is mean-square stable, then the modified equation \cref{eq:smrtesteqmod} is mean-square stable for any value of $\eta\geq 0$.
\end{theorem}
\begin{proof}
The modified equation \cref{eq:smrtesteqmod} is mean-square stable if, and only if, 
\begin{equation}\label{eq:stabcond_smrtesteqmod}
	\varphi(\eta\lambda)(\lambda+\zeta)+\frac{1}{2}|\varphi(\eta\lambda/2)\mu|^2< 0.
\end{equation}
Using \cref{lemma:phisq} and $(\lambda,\zeta,\mu)\in \mathcal{S}^{mMS}$ it follows
\begin{equation}
\varphi(\eta\lambda)(\lambda+\zeta)+\frac{1}{2}|\varphi(\eta\lambda/2)\mu|^2\leq \varphi(\eta\lambda)(\lambda+\zeta+\frac{1}{2}|\mu|^2)< 0. \tag*{\qedhere}
\end{equation}
\end{proof}
In view of \cref{lemma:phig,lemma:phisq,thm:msstab_modsde} we understand why we need a factor $1/2$ in the definition of $\ge$ in \cref{eq:defv}; this guarantees the right damping for the diffusion term. In practice, $\eta$ is chosen so that $\eta|\zeta|\geq 2$ as from \cref{thm:eta} this choice of $\eta$ guarantees that stiffness of \cref{eq:modsde} depends only on the slow term $\fs$. In \cref{fig:phi} we illustrate the inequality $\varphi(z/2)^2\leq \varphi(z)$ and see that it is very tight; hence, replacing $g$ by $\ge$ guarantees mean-square stability without over damping the diffusion term.
In \cref{fig:origmsde,fig:modsdea} we illustrate the stability conditions $(\lambda,\zeta,\mu)\in \mathcal{S}^{mMS}$ and \cref{eq:stabcond_smrtesteqmod}, respectively, and show that if $(\lambda,\zeta,\mu)\in \mathcal{S}^{mMS}$ is satisfied then also \cref{eq:stabcond_smrtesteqmod} is satisfied. We however emphasize, as illustrated in \cref{fig:modsdeb}, that if the noise term is not damped the modified equation might be unstable.

\section{The multirate second-kind orthogonal Runge--Kutta--Chebyshev method}\label{sec:mskrock}
We introduce here a stabilized explicit multirate method for \cref{eq:msde} based on the stochastic modified equation \cref{eq:modsde}: the mSK-ROCK scheme. We first recall the mRKC method for the deterministic multirate differential equation \cref{eq:mrode} based on the modified equation \cref{eq:modode}.

\subsection{The multirate Runge--Kutta--Chebyshev method}\label{sec:mrkc}
The multirate Runge--Kutta--Chebyshev (mRKC) method is obtained by discretizing \eqref{eq:modode} with an $s$-stage Runge--Kutta--Chebyshev (RKC) method and approximating $\fe$ given in \cref{def:fe} by solving \eqref{eq:defu} with an $m$-stage RKC method. The RKC method \cite{HoS80} employed for the approximation of \cref{eq:modode,eq:defu} is a stabilized explicit scheme with stability domain growing quadratically with the number of function evaluations, see \cite{SSV98,HoS80,VHS90} for more details.

\subsubsection*{The algorithm}
Let $\tau>0$ be the step size and the stages $s,m$ of the two RKC methods satisfy the stability conditions
\begin{align}\label{eq:defsmeta}
\tau\rhos &\leq  \beta s^2, & \eta\rhof &\leq  \beta m^2, & \mbox{with}&& \eta&=\frac{6\tau}{\beta s^2}\frac{m^2}{m^2-1},
\end{align}
$\beta=2-4/3\varepsilon$ and typically $\varepsilon=0.05$. The value of $\eta$ follows from the stability analysis of the scheme, see \cite{AGR20}. One step of the $\mRKC$ scheme is given by a classical RKC scheme applied to the modified equation \cref{eq:modode}, i.e.,
\begin{equation}\label{eq:defmrkc}
\begin{aligned}
k_0&=y_{n},\\
k_1 &= k_0+\mu_1\tau\bfe(k_0),\\
k_j&= \nu_j k_{j-1}+\kappa_j k_{j-2}+\mu_j\tau \bfe(k_{j-1}) \quad j=2,\ldots,s,\\
y_{n+1} &=k_s,
\end{aligned}
\end{equation}
where the coefficients depend on $\oz=1+\varepsilon/s^2$,  $\ou=T_s(\oz)/T_s'(\oz)$, with $T_s(x)$ the Chebyshev polynomial of the first kind of degree $s$. 
The parameters of the scheme are given by $\mu_1 = \ou/\oz$,
\begin{align}\label{eq:defcoeff} 
	\mu_j&= 2\ou  b_j/b_{j-1}, & 
	\nu_j&= 2\oz b_j/b_{j-1}, & 
	\kappa_j&=-b_j/b_{j-2},  &
	\text{for }j&=2,\ldots,s,
\end{align}
with $b_j=T_j(\oz)^{-1}$ for $j=0,\ldots,s$.
\modr{Note that in \cref{eq:defmrkc} only three vectors must be stored (even if $s$ is large). The} recurrence relation allows also for a good internal stability with respect to roundoff errors \cite{VHS90}. Moreover, the scheme \cref{eq:defmrkc} is stable for $\tau\rhos\leq\beta s^2$ and this ensures a quadratic growth of the stability domain with respect to the stage number $s$. This is in sharp contrast with the explicit Euler method (where the stability domain grows only linearly with respect to the number of steps).

Following \eqref{eq:deffedif} the averaged force is defined by
\begin{equation}\label{eq:defbfe}
\bfe(y)=\frac{1}{\eta}(u_\eta-y),
\end{equation}
where $u_\eta$ is obtained applying one step of size $\eta$ of a RKC method with $m$ stages to \eqref{eq:defu}. Hence, it is computed with the scheme
\begin{equation}\label{eq:defbue}
\begin{aligned}
u_0&= y, \qquad  u_1 = u_0+\alpha_1\eta(\ff(u_0)+\fs(y)),\\
u_j &= \beta_j u_{j-1}+\gamma_j u_{j-2}+\alpha_j\eta (\ff(u_{j-1})+\fs(y)) \quad j=2,\ldots,m,\\
u_\eta &=u_m.
\end{aligned}
\end{equation}
The parameters of the $m$-stage RKC scheme are given by $\vz=1+\varepsilon/m^2$, $\vu=T_m(\vz)/T_m'(\vz)$, $a_j=T_j(\vz)^{-1}$ for $j=0,\ldots,m$ and 
%\begin{align}\label{eq:defcoeff2} 
%\alpha_1 &= \vu/\vz, & \alpha_j&= 2\vu  a_j/a_{j-1}, \\ \notag
%\beta_j&= 2\vz a_j/a_{j-1}, & \gamma_j&=-a_j/a_{j-2},  &\text{for }j&=2,\ldots,m.
%\end{align}
$\alpha_1 = \vu/\vz$, 
\begin{align}\label{eq:defcoeff2} 
	\alpha_j&= 2\vu  a_j/a_{j-1}, &
	\beta_j&= 2\vz a_j/a_{j-1}, & 
	\gamma_j&=-a_j/a_{j-2},  &
	\text{for }j&=2,\ldots,m.
\end{align}
The $\mRKC$ scheme \cref{eq:defsmeta,eq:defmrkc,eq:defcoeff,eq:defbfe,eq:defbue,eq:defcoeff2} is first-order accurate \cite{AGR20}.

\subsubsection*{Linear stability analysis on the multirate test equation}
Here we recall the stability properties of the $\mRKC$ scheme, as they are crucial for studying the stability of the mSK-ROCK scheme introduced in \cref{sec:mskrockalg}. First, we compute a closed expression for $\bfe$ when the $\mRKC$ scheme \cref{eq:defsmeta,eq:defmrkc,eq:defcoeff,eq:defbfe,eq:defbue,eq:defcoeff2} is applied to the multirate test equation \cref{eq:mrtesteq}. Let
\begin{equation}\label{eq:defPhim}
A_m(z)=\frac{T_m(\vz+\vu z)}{T_m(\vz)},\qquad\qquad \Phi_m(z)=\frac{A_m(z)-1}{z} \quad\mbox{for }z\neq 0
\end{equation}
and $\Phi_m(0)=1$, where $A_m(z)$ is the stability polynomial of the $m$-stage RKC scheme satisfying $|A_m(z)|\leq 1$ for $|z|\leq\beta m^2$. The function $\Phi_m(z)$ is the numerical counterpart of $\varphi(z)$ given in \cref{eq:defphi}, indeed in \modr{\cite[Section 4.1]{AGR20}} it is proved the following.
\begin{lemma}\label{lemma:closedbue}
	Let $\lambda,\zeta\leq 0$, $\fs(y)=\zeta y$, $\ff(y)=\lambda y$, $\eta>0$, $m\in\Nb$ and $y\in\Rb$. Then
	\begin{equation}\label{eq:bfetesteq}
	\bfe(y) =\Ps_m(\eta\lambda)(\lambda+\zeta)y.
	\end{equation}
\end{lemma}
Therefore, as $\varphi(z)$, $\Phi_m(z)$ has a smoothing effect on $\lambda+\zeta$, which decreases the stiffness of the problem as long as $|z|\leq \beta m^2$ and thus $|A_m(z)|\leq 1$.
Plugging $\bfe(y)$ from \cref {eq:bfetesteq} into \cref {eq:defmrkc} leads to 
\begin{equation}
y_{n+1} = A_s(\tau \Phi_m(\eta\lambda)(\lambda+\zeta))y_{n},
\end{equation} 
where $A_s(p)$ is the stability polynomial of the $s$-stage RKC scheme. Hence, the scheme is stable if $|\tau \Phi_m(\eta\lambda)(\lambda+\zeta)|\leq \beta s^2$ and thus $|A_s(\tau \Phi_m(\eta\lambda)(\lambda+\zeta))|\leq 1$. In \modr{\cite[Theorem 4.5]{AGR20}} the following result is proved.
\begin{theorem}\label{thm:stab_mrkc}
	Let the damping $\varepsilon=0$, $\lambda\leq 0$ and $\zeta<0$. Then, for all $\tau>0,s,m$ and $\eta$ satisfying \cref{eq:defsmeta} with $\rhof=|\lambda|$ and $\rhos=|\zeta|$
%	\begin{align}\label{eq:defsmetascalar}
%	\tau |\zeta| &\leq\beta s^2, & \eta|\lambda| &\leq \beta m^2 &\text{with}  && \eta &\geq \frac{6\tau}{\beta s^2}\frac{m^2}{m^2-1},
%	\end{align}
	it holds $|A_s(\tau \Phi_m(\eta\lambda)(\lambda+\zeta))|\leq 1$, i.e. the $\mRKC$ scheme is stable.
\end{theorem}
It is shown in \cite{AGR20} that the scheme is stable also for small damping parameters $\varepsilon>0$ and numerical experiments confirm that stability holds for any damping, here we consider $\varepsilon=0$ for simplicity.

\subsection{The multirate SK-ROCK method}\label{sec:mskrockalg}
The $\mSKROCK$ method is a generalization of the mRKC method of \cref{sec:mrkc} to SDEs. It consists in the time discretization of \eqref{eq:modsde} with the SK-ROCK scheme \cite{AAV18}, but where $\fe$ is replaced by $\bfe$ and $\ge$ of \cref{def:ge} is approximated by solving the two auxiliary problems in \cref{eq:defv} with a modified RKC method.

\subsubsection*{The algorithm}
For simplicity, we define the mSK-ROCK method for a vector valued diffusion term $g:\Rn\rightarrow\Rn$ and generalize the scheme to a matrix valued diffusion $g:\Rn\rightarrow\Rb^{n\times l}$ at the end of this section.

Let $s$, $m$ and $\eta$ be as in \eqref{eq:defsmeta} but with the constraint that $m$ must be even. Denote $m=2r$ with \modr{$r\in\Nb\setminus\{0\}$}. One step of the $\mSKROCK$ method is given by
\begin{equation}\label{eq:mskrock}
\begin{aligned}
K_0&=X_n,\\
K_1 &= K_0+\mu_1\tau \bfe(K_0+\nu_1\bQe )+\kappa_1 \bQe,\\
K_j&= \nu_j K_{j-1}+\kappa_j K_{j-2}+\mu_j\tau \bfe(K_{j-1}) \quad j=2,\ldots,s,\\
X_{n+1}&=K_s,
\end{aligned}
\end{equation}
with $\bQe = \bge(K_0)\Delta W_n$, $\Delta W_n=W(\tnpu)-W(\tn)$, $\bfe$ as in \cref{eq:defbfe},\cref{eq:defbue}, $\mu_j,\nu_j,\kappa_j$ for $j=2,\ldots,s$ as in \eqref{eq:defcoeff} and 
$\mu_1=\ou/\oz$, $\nu_1=s\,\ou/2$, $\kappa_1= s\,\ou/\oz$. Observe that in this scheme the noise term is introduced in the first stage. 

The function $\bge$ is a numerical approximation of $\ge$. From \eqref{eq:defgedif} we define
\begin{equation}\label{eq:defbge}
\bge(x)=\frac{1}{\eta}(v_\eta-\bv _\eta),
\end{equation}
where $v_\eta$ and $\bv _\eta$ are approximations of $v(\eta)$ and $\bv (\eta)$, respectively. We compute $v_\eta$ using a modified $r$-stage RKC scheme: the parameters are those of a RKC scheme with $m=2r$ stages and the contribution of $g(x)$ appears only in the first stage. Hence, $v_\eta$ is given by
\begin{equation}\label{eq:defve}
\begin{aligned}
v_0 &= x,\qquad 
v_1 = v_0+\alpha_1\eta \ff(v_0+\beta_1 \theta_1 \eta g(x))+\gamma_1\theta_1 \eta g(x),\\
v_j &= \beta_j v_{j-1}+\gamma_j v_{j-2}+\alpha_j\eta \ff(v_{j-1}) \quad j=2,\ldots,r,\\
v_\eta &=v_r
\end{aligned}
\end{equation}
and $\bv _\eta$ is given by
\begin{equation}\label{eq:defbve}
\begin{aligned}
\bv_0 &= x,\qquad 
\bv _1 = \bv _0+\alpha_1\eta \ff(\bv _0),\\
\bv _j &= \beta_j \bv _{j-1}+\gamma_j \bv _{j-2}+\alpha_j\eta \ff(\bv _{j-1}) \quad j=2,\ldots,r,\\
\bv _\eta &=\bv _r,
\end{aligned}
\end{equation}
where in \cref{eq:defve,eq:defbve} the parameters $\alpha_1$ and $\alpha_j,\beta_j,\gamma_j$ for $j=2,\ldots,r$ are the parameters of the $m$-stage RKC scheme given in \eqref{eq:defcoeff2} with $m=2r$ and the additional parameters in \cref{eq:defve} are given by $\beta_1=m\vu/2$, $\gamma_1= m\,\vu/\vz$ and $\theta_1 = T_r(\vz)/(2\vu T_r'(\vz))$.
%\begin{equation}\label{eq:deftheta}
%\beta_1=m\frac{\vu}{2},\qquad\qquad \gamma_1= m\frac{\vu}{\vz}, \qquad\qquad \theta_1 = \frac{1}{2\vu}\frac{T_r(\vz)}{T_r'(\vz)}.
%\end{equation}

Note that the $1/2$ factor in \eqref{eq:defv} disappears from \labelcref{eq:defve,eq:defbve} but is reflected on the fact that we take $r=m/2$ stages. Indeed, for the approximation $\bge$ of $\ge$ for linear problems (see \cref{eq:defgephi}) the numerical counterpart of \cref{lemma:phisq} rely on the identity $2T_r(x)^2=T_{2r}(x)+1$ of Chebyshev polynomials. This sets the relation between the number of stages $r$ for computing $\bge$ and and the number of stages $m=2r$ for computing $\bfe$, see \cref{lemma:relBPhi} below.

Now, we discuss the case where $g:\Rn\rightarrow\Rb^{n\times l} $ is a matrix valued function and therefore \cref{eq:defve} is not well-defined. One possible approach is to compute \cref{eq:defve} for each column of $g$ and build a modified matrix $\ge$ column-wise. 
However, this way of proceeding entails the computation of \cref{eq:defve} for each column of $g$, which can rapidly become expensive. A better solution is to replace $\eta g(x)$ in \cref{eq:defve} by $\eta g(x)\Delta W_n$, which is vector valued and therefore \cref{eq:defbge,eq:defve,eq:defbve} can be computed. Then we replace $\bQe$ in \cref{eq:mskrock} by $\bQe=\bge(x)$, as $\Delta W_n$ is already contained in $\bge(x)$. With this second approach, \cref{eq:defve} is computed only once and therefore the cost of stabilizing a vector or a matrix valued diffusion term $g$ is equivalent. Note that when $\ff$ is linear the two approaches give exactly the same result $\bge$. When $\ff$ is nonlinear we obtain two slightly different methods, nonetheless we can show that both have the same accuracy and mean-square stability properties.

\subsubsection*{Efficiency analysis}
Given the spectral radii $\rhof,\rhos$ of the Jacobians of $\ff,\fs$, respectively, we want to compare the theoretical efficiency, in terms of function evaluations, of the $\mSKROCK$ and SK-ROCK method.
We set $\varepsilon=0$ and let $s,m$ vary in $\Rb$. The cost of evaluating $\ff,\fs,g$ relatively to the cost of evaluating $\ff+\fs+g$ is denoted $c_F,c_S,c_g\in [0,1]$, respectively, with $c_F+c_S+c_g=1$. 

One step of $\mSKROCK$ requires $s$ evaluations of $\bfe$ and one of $\bge$. Each evaluation of $\bfe$ needs $m$ evaluations of $\ff$ and one of $\fs$, an evaluation of $\bge$ requires $2r=m$ evaluations of $\ff$ and one of $g$. Hence, the cost of one step of $\mSKROCK$ is given by 
\begin{equation}
C_{\mSKROCKop}= s(m c_F+c_S)+m c_F+c_g=((s+1)m-1)c_F+(s-1)c_S+1,
\end{equation}
where we used $c_g=1-c_F-c_S$.
Conditions \cref{eq:defsmeta} with $\beta=2$ yield $s=\sqrt{\tau\rhos/2}$ and $m=\sqrt{3\rhof/\rhos+1}$, thus
\begin{equation}
C_{\mSKROCKop}=\left(\left(\sqrt{\frac{\tau\rhos}{2}}+1\right)\sqrt{3\frac{\rhof}{\rhos}+1}-1\right)c_F+\left(\sqrt{\frac{\tau\rhos}{2}}-1\right)c_S+1.
\end{equation}
In contrast, the standard SK-ROCK method is given by \cref{eq:mskrock} but with $\bfe$ replaced by $f=\ff+\fs$ and $\bge$ replaced by $g$. Hence, one step of SK-ROCK needs $s$ evaluations of $\ff+\fs$ and one of $g$, with $s=\sqrt{\tau\rho/2}$, where $\rho$ is the spectral radius of the Jacobian of $f$ and we assume $\rho=\rhof+\rhos$. Thus, the cost of one step of SK-ROCK is
\begin{equation}
	C_{\SKROCKop}= s(c_F+c_S)+c_g=(s-1)(c_F+c_S)+1 = \left(\sqrt{\frac{\tau(\rhof+\rhos)}{2}}-1\right)(c_F+c_S)+1.
\end{equation}

\begin{figure}%[!tbp]
\begin{center}
\begin{subfigure}[t]{\subfigsize\textwidth}
\centering
\begin{tikzpicture}
\node at (0,0) {\includegraphics[trim=0cm 0cm 0cm 0cm, clip, width=0.8\textwidth]{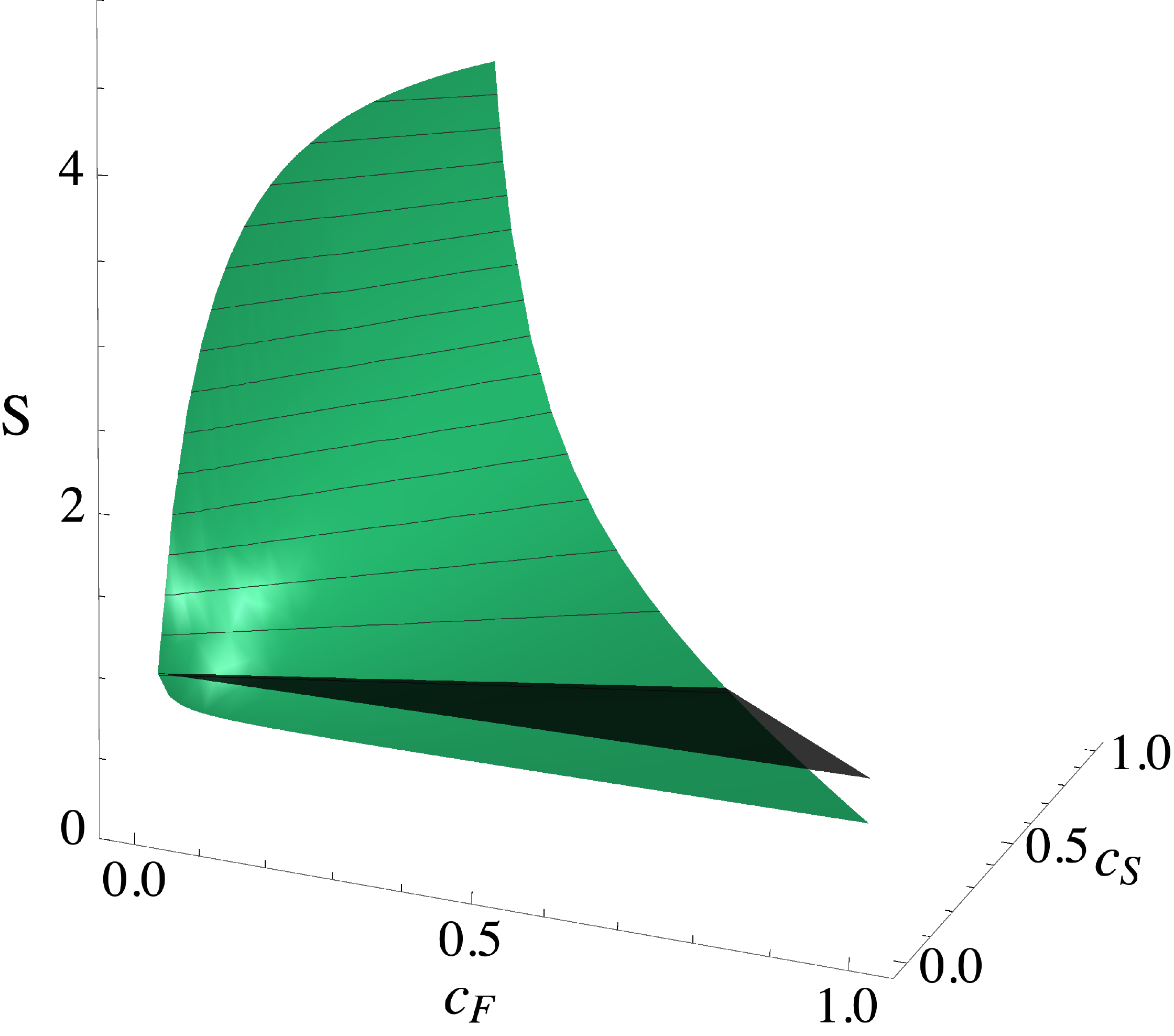}};
\draw[colorone, fill] (0.8,1) circle (0.1cm) node[anchor=west,black] {$\;S$};
\draw[black, fill,opacity=0.7] (0.8,0.5) circle (0.1cm) node[anchor=west,opacity=1] {$\;1$};
\end{tikzpicture}
\caption{Speed-up for $p_F=2000$, $p_S=200$.}
\label{fig:speedsrkcsqa}
\end{subfigure}
\hfill
\begin{subfigure}[t]{\subfigsize\textwidth}
\centering
\begin{tikzpicture}
\node at (0,0) {\includegraphics[trim=0cm 0cm 0cm 0cm, clip, width=0.8\textwidth]{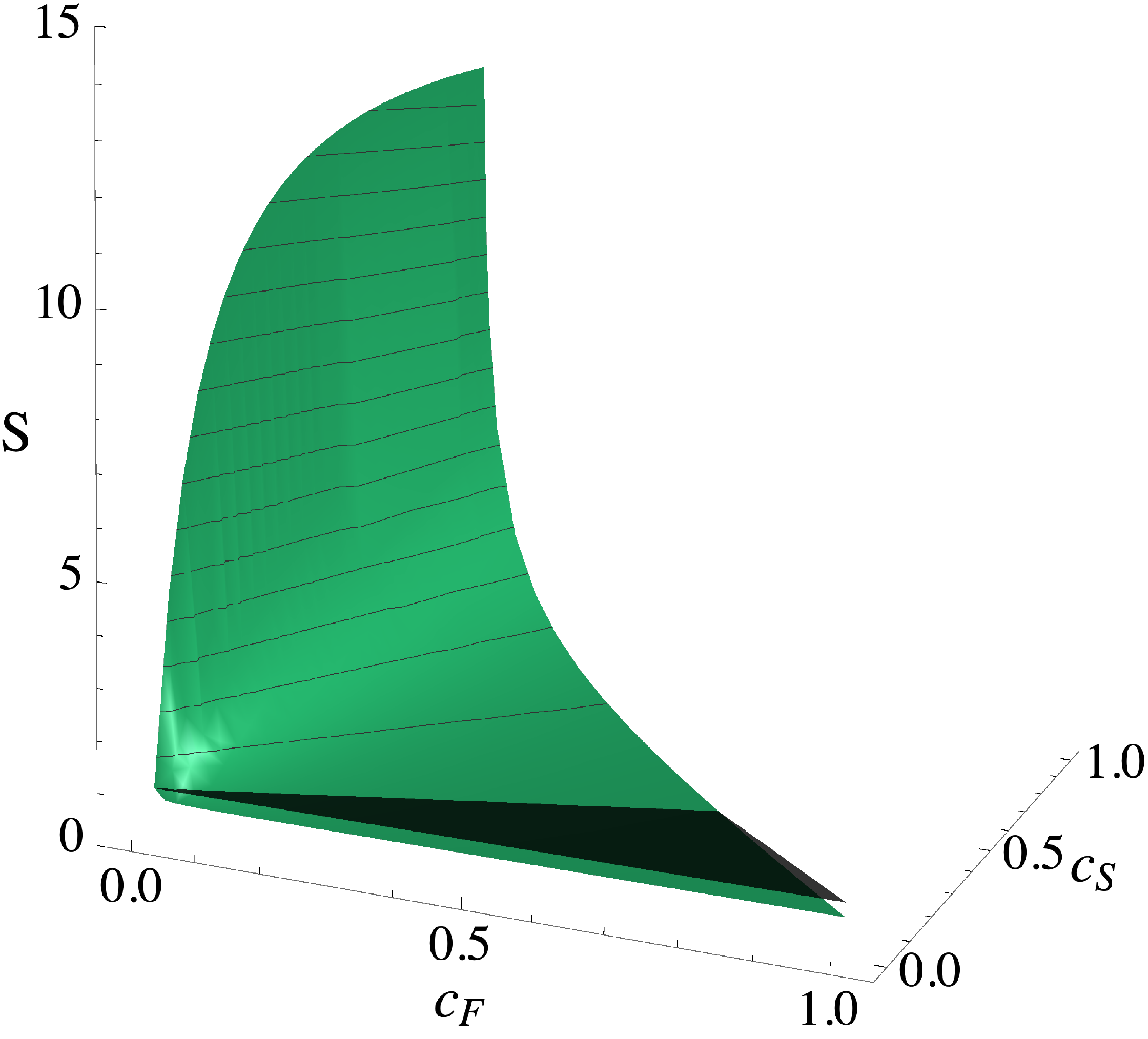}};
\draw[colorone, fill] (0.8,1) circle (0.1cm) node[anchor=west,black] {$\;S$};
\draw[black, fill,opacity=0.7] (0.8,0.5) circle (0.1cm) node[anchor=west,opacity=1] {$\;1$};
\end{tikzpicture}
\caption{Speed-up for $p_F=20000$, $p_S=200$.}
\label{fig:speedsrkcsqb}
\end{subfigure}
\end{center}
\caption{The relative speed-up $S$ of $\mSKROCK$ over SK-ROCK, for some fixed values of $p_F=\tau\rhof$ and $p_S=\tau\rhos$. The black surface corresponds to the constant function $1$, for $S>1$ the $\mSKROCK$ scheme is faster than SK-ROCK while for $S<1$ the SK-ROCK scheme is faster than $\mSKROCK$.}
\label{fig:speedsrkcsq}
\end{figure}
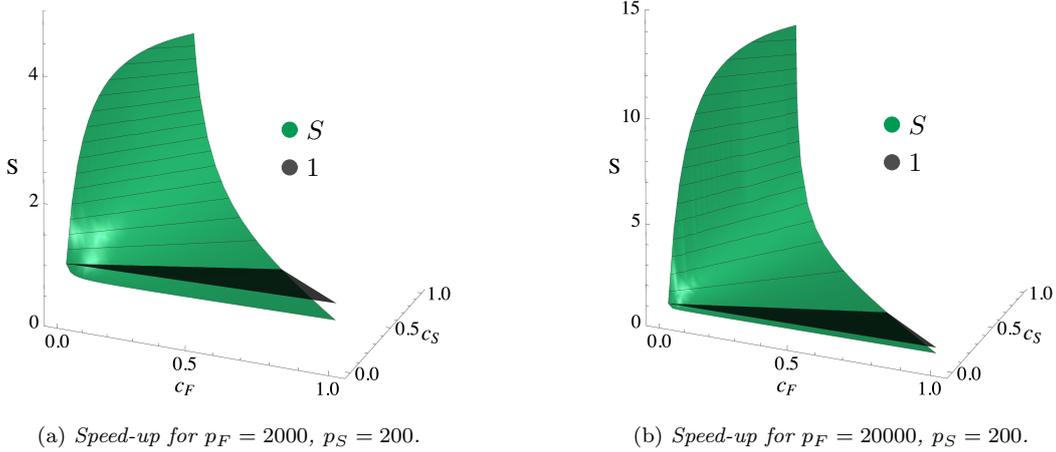

Let $p_F=\tau\rhof$ and $p_S=\tau\rhos$, the theoretical relative speed-up $S$ is defined as the ratio between the two costs:
\begin{equation}
S = \frac{C_{\SKROCKop}}{C_{\mSKROCKop}}=\frac{\left(\sqrt{p_F+p_S}-\sqrt{2}\right)(c_F+c_S)+\sqrt{2}}{\left(\left(\sqrt{p_S}+\sqrt{2}\right)\sqrt{3\frac{p_F}{p_S}+1}-\sqrt{2}\right)c_F+\left(\sqrt{p_S}-\sqrt{2}\right)c_S+\sqrt{2}}.
\end{equation}
For some values of $p_F,p_S$ we display $S$ in function of $c_F$, $c_S$ in \cref{fig:speedsrkcsq}, with $c_F+c_S\in[0,1]$. 
We see that the speed-up increases as $c_S\to 1$; indeed, SK-ROCK needs more evaluations of $\fs$. In contrast, the $\mSKROCK$ method is slower than SK-ROCK ($S<1$) if $c_F$ is not sufficiently small, as it needs more evaluations of $\ff$. However, we recall that we intend to use the $\mSKROCK$ method when $\ff$ is cheap to evaluate, otherwise we simply use the SK-ROCK scheme. 

\section{Stability and convergence analysis}\label{sec:analysis}
This section is devoted to the stability and accuracy analysis of the $\mSKROCK$ method.
%, given by \cref{eq:mskrock,eq:defbge,eq:defve,eq:defbve,eq:defbfe,eq:defbue}.

\subsection{Stability analysis on the stochastic multirate test equation}\label{sec:stabanalyis}
We show here that when the $\mSKROCK$ method is applied to the stochastic multirate test equation \eqref{eq:smrtesteq} the scheme is stable. \modr{The stability analysis presented here, on the stochastic multirate test equation, is straightforwardly generalizable only to problems for which the Jacobians of $\ff,\fs,g$ are simultaneously diagonalizable. Otherwise, multidimensional noise can have effects that are not captured by a scalar dynamics and model problems in higher dimension must be considered, as in \cite{BuK10}. However, the results of \cite{BuK10} also suggest that the scalar model problem offers a first practical insight on the stability of more general problems. Hence, here we consider only the scalar linear problem. Also, we believe that doing the stability analysis on a more complex model is out of the scope of this paper.}

In order to analyze the stability of the $\mRKC$ method in \cref{sec:mrkc} we computed a closed expression for $\bfe$ in \modr{\cref{lemma:closedbue}}. We start by deriving an expression for $\bge$ given in the next lemma. Define 
\begin{equation}\label{eq:defPsir}
\Psi_r(z)= \frac{U_{r-1}(\vz+\vu z)}{U_{r-1}(\vz)}\left(1+\frac{\vu}{2}z\right),
\end{equation}
where $U_k(x)$ is the Chebyshev polynomial of the second kind of degree $k$ and $\vz,\vu$ are given in \cref{sec:mrkc}. The Chebyshev polynomials of the second kind have a recurrence relation $U_k(x)=2x U_{k-1}(x)-U_{k-2}(x)$, similar to the Chebyshev polynomials of the first kind $T_k(x)$ except for the initial values $U_0(x)=1$ and $U_1(x)=2x$.
\begin{lemma}\label{lemma:closedbve}
Under the assumptions of \cref{lemma:closedbue} and $\ge(x)=\mu x$ with $\mu,x\in \Rb$, it holds
\begin{equation}\label{eq:bgetesteq}
\bge(x) =\Psi_r(\eta\lambda)\mu x.
\end{equation}
\end{lemma}
\begin{proof} %shorter, long verson below
	Replacing $\ff(x)=\lambda x$ and $g(x)=\mu x$ in \eqref{eq:defve} yields
	\begin{align}\label{eq:vetesteqa}
		\begin{split}
			v_0&=x,\qquad
			v_1 = v_0+\alpha_1 z v_0+ r_1,\\
			v_j &= \beta_j v_{j-1}+\gamma_j v_{j-2}+\alpha_j z v_{j-1} \quad j=2,\ldots,r,
		\end{split}
	\end{align}
	with $z=\eta\lambda$ and $r_1=(\alpha_1\beta_1 z +\gamma_1)\theta_1 \eta \mu x$. Scheme \eqref{eq:vetesteqa} is a perturbed (by $r_1$) RKC scheme, these schemes are studied in \cite{VHS90} and it can be shown (see \cite{Ros20} for more details) that
	\begin{equation}\label{eq:vetesteqb}
			v_\eta = a_r T_r(\vz+\vu z) x + \frac{a_r}{a_1}U_{r-1}(\vz+\vu z)r_1.
	\end{equation}
	From the definitions of $\alpha_1,\beta_1,\gamma_1,\theta_1,\vz,a_j$ and $T_1(\vz)=\vz$, $T_r'(\vz)=rU_{r-1}(\vz)$ it follows
	\begin{align}
		\frac{a_r}{a_1}r_1 &= \frac{1}{U_{r-1}(\vz)} \left(1+\frac{\vu}{2}z \right)\eta\mu x
	\end{align}
	and thus $v_\eta = a_r T_r(\vz+\vu z) x +\eta \Psi_r(z)\mu x$.
	We have as well $\bv _\eta = a_r T_r(\vz+\vu z) x$, which, with \cref{eq:defbge}, yields \cref{eq:bgetesteq}.
\end{proof}
We now apply the $\mSKROCK$ method to the stochastic multirate test equation \cref{eq:smrtesteq}. Let $\Delta W_n=\tau^{1/2}\xi_n$ with $\xi_n\sim \mathcal{N}(0,1)$, plugging \cref{eq:bfetesteq,eq:bgetesteq} into \cref{eq:mskrock} yields
\begin{equation}
	X_{n+1}= R_s(p_m,q_r,\xi)X_n, \quad\mbox{where}\quad R_s(p,q,\xi)=A_s(p)+B_s(p)q\xi
\end{equation}
is the stability polynomial of the SK-ROCK method \cite{AAV18}, with $A_s(p)$ given by \cref{eq:defPhim} and
\begin{equation}
	B_s(p)=\frac{U_{s-1}(\oz+\ou p)}{U_{s-1}(\oz)}\left(1+\frac{\ou}{2}p\right),
	\quad p_m=\tau\Ps_m(\eta\lambda)(\lambda+\zeta), \quad q_r= \Psi_r(\eta\lambda)\mu\tau^{1/2}.
\end{equation}
%and
%\begin{equation}\label{eq:pmqr}
%p_m=\tau\Ps_m(\eta\lambda)(\lambda+\zeta),\qquad\mbox{and}\qquad q_r= \Psi_r(\eta\lambda)\mu\tau^{1/2}.
%\end{equation}
%This motivates the following definition.
%\begin{definition}\label{def:stabpolmskrock}
%Let $\tau>0$, $\eta\geq 0$, $\lambda,\zeta\leq 0$, $\mu\in\Rn$ and $\xi\sim\mathcal{N}(0,1)$ a Gaussian random variable. The stability polynomial of the $(s,m)$-stage $\mSKROCK$ method is defined as
%\begin{equation}
%R_{s,m}(\tau,\eta,\lambda,\zeta,\mu,\xi)=A_s(p_m)+B_s(p_m)q_r\xi,
%\end{equation}
%with $A_s$, $B_s$ as in \eqref{eq:defAB}, $p_m,q_r$ as in \eqref{eq:pmqr} and $2r=m$.
%% and $\Ps_m$, $\Psi_r$ as in \cref{eq:defPhim,eq:defPsir}, respectively.
%\end{definition}
The next lemma is the numerical counterpart of \cref{lemma:phisq} and therefore it is the main tool for proving stability of the scheme in \cref{thm:stab_mskrock} below.
\begin{lemma}\label{lemma:relBPhi}
Let $m,r\in\Nb^*$, $m=2r$ and $\varepsilon=0$. Then $\Psi_r(z)^2\leq \Ps_m(z)$ for all $z\in [-\beta m^2,0]$.
\end{lemma}
\begin{proof}
Since $\Psi_r(0)^2=\Phi_m(0)=1$ we consider $z\in [-\beta m^2,0)$, with $\beta=2$. For $\varepsilon=0$ we also have $\vz=1$, $\vu=1/m^2$. Letting $x=\vz+\vu z=1+z/m^2\in [-1,1)$ and using $U_{r-1}(1)=r$ yields
\begin{equation}\label{eq:psix}
\Psi_r(z)=\frac{U_{r-1}(x)}{2r}(x+1).
\end{equation}
The identity $2T_k(x)T_j(x)=T_{k+j}(x)+T_{|k-j|}(x)$ implies $T_{2r}(x)=2T_r(x)^2-1$ and thus
\begin{equation}\label{eq:phix}
\Phi_{2r}(z)=\frac{T_{2r}(x)-1}{z}=\frac{2T_{r}(x)^2-2}{(2r)^2(x-1)}=\frac{T_r(x)^2-1}{2r^2(x-1)}.
\end{equation}
From \cref{eq:psix}, \cref{eq:phix} and $x-1<0$, $\Psi_r(z)^2\leq \Ps_{2r} (z)$ is equivalent to
\begin{equation}
0\leq U_{r-1}(x)^2(x^2-1)(x+1)-2(T_r(x)^2-1).
\end{equation}
The result follows from the identity $T_r(x)^2-1=U_{r-1}(x)^2(x^2-1)$ and $x\in [-1,1)$.
\end{proof}
Numerical evidences show that \cref{lemma:relBPhi} is valid for any damping parameter $\varepsilon\geq 0$. Indeed, we display $\Ps_{2r}(z)$ and $\Psi_r(z)^2$ for $r=3$ in \cref{fig:phipsisrkcsq}, for a small damping $\varepsilon=0.05$ and a high damping $\varepsilon=1$. In both cases relation $\Psi_r(z)^2\leq\Phi_{2r}(z)$ holds and is tight. 
\begin{figure}
	\begin{center}
		\begin{subfigure}[t]{\subfigsize\textwidth}
			\centering
			\begin{tikzpicture}[scale=\plotimscale]
			\begin{axis}[height=\aspectratio*\plotimsized\textwidth,width=\plotimsized\textwidth, ymin=-0.2, ymax=1,xmax=0,xmin=-72,legend columns=1,legend style={draw=\legendboxdraw,fill=\legendboxfill,at={(0,1)},anchor=north west},legend cell align={left},
			xlabel={$z$}, ylabel={},label style={font=\normalsize},tick label style={font=\normalsize},legend image post style={scale=\legendmarkscale},legend style={nodes={scale=\legendfontscale, transform shape}},grid=none]
			\addplot[color=colorone,line width=\plotlinewidth pt,mark=\markone,mark repeat=20,mark phase=0,mark size=\plotmarksizeu pt] table [x=z,y=phi,col sep=comma] 
			{data/text/phi_psisq_r_3.csv};\addlegendentry{$\Ps_6(z)$}
			\addplot[color=colortwo,line width=\plotlinewidth pt,mark=\marktwo,mark repeat=20,mark phase=10,mark size=\plotmarksizeu pt] table [x=z,y=psisq,col sep=comma] 
			{data/text/phi_psisq_r_3.csv};\addlegendentry{$\Psi_3(z)^2$}
			\coordinate (pta) at (axis cs:-17,0.055);
			\coordinate (ptb) at (axis cs:-23,0.25);
			\end{axis}
			\draw[] (pta)node[circle,anchor=north west,draw]{}--(ptb)node[anchor=south east,draw,circle,fill=white,inner sep=0,outer sep=0]{
				\begin{tikzpicture}[scale=\zoomsize,trim axis left,trim axis right]
				\begin{axis}[tiny,hide axis,xlabel={},ylabel={},ticks=none,xmin=-22,xmax=-14,ymin=-0.005,ymax=0.02]
				\addplot[color=colorone,line width=\zoomlinewidth pt,mark=\markone,mark repeat=10,mark phase=0,mark size=4 pt] table [x=z,y=phi,col sep=comma] 
				{data/text/phi_psisq_r_3.csv};
				\addplot[color=colortwo,line width=\zoomlinewidth pt,mark=\marktwo,mark repeat=10,mark phase=5,mark size=4 pt] table [x=z,y=psisq,col sep=comma] 
				{data/text/phi_psisq_r_3.csv};
				\end{axis}
				\end{tikzpicture}
			};
		\end{tikzpicture}
		\caption{Small damping $\varepsilon=0.05$.}
	\end{subfigure}\hfill
	\begin{subfigure}[t]{\subfigsize\textwidth}
		\centering
		\begin{tikzpicture}[scale=\plotimscale]
		\begin{axis}[height=\aspectratio*\plotimsized\textwidth,width=\plotimsized\textwidth, ymin=-0.2, ymax=1,xmax=0,xmin=-45.5,legend columns=1,legend style={draw=\legendboxdraw,fill=\legendboxfill,at={(0,1)},anchor=north west},legend cell align={left},
		xlabel={$z$}, ylabel={},label style={font=\normalsize},tick label style={font=\normalsize},legend image post style={scale=\legendmarkscale},legend style={nodes={scale=\legendfontscale, transform shape}},grid=none]
		\addplot[color=colorone,line width=\plotlinewidth pt,mark=\markone,mark repeat=20,mark phase=0,mark size=\plotmarksizeu pt] table [x=z,y=phi,col sep=comma] 
		{data/text/phi_psisq_r_3_high_damp.csv};\addlegendentry{$\Ps_{6}(z)$}
		\addplot[color=colortwo,line width=\plotlinewidth pt,mark=\marktwo,mark repeat=20,mark phase=10,mark size=\plotmarksizeu pt] table [x=z,y=psisq,col sep=comma] 
		{data/text/phi_psisq_r_3_high_damp.csv};\addlegendentry{$\Psi_3(z)^2$}
		\end{axis}
		\end{tikzpicture}
		\caption{High damping $\varepsilon=1$.}
	\end{subfigure}
\end{center}
\caption{Illustration of $\Ps_{2r}(z)$ and $\Psi_r(z)^2$ for $r=3$, with damping $\varepsilon=0.05$ (left) or $\varepsilon=1$ (right).}
\label{fig:phipsisrkcsq}
\end{figure}
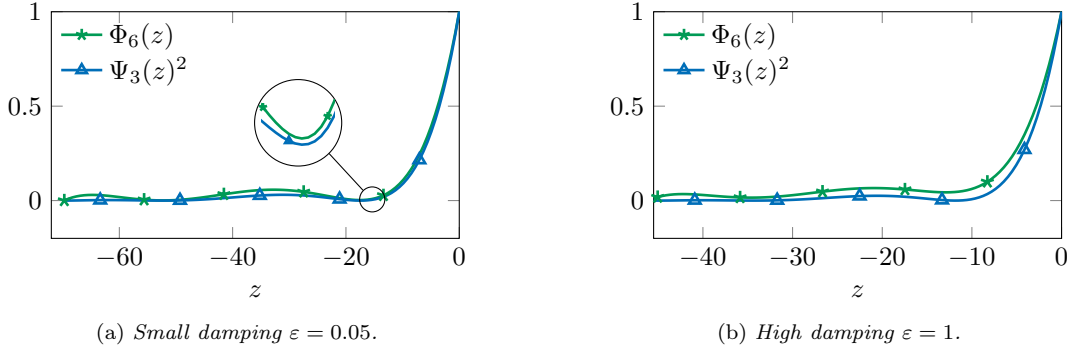
\begin{theorem}\label{thm:stab_mskrock}
Let $\varepsilon=0$ and $(\lambda,\zeta,\mu)\in\mathcal{S}^{mMS}$ (see \cref{eq:defSmMS}).
Then, for all $\tau>0,s,m$ and $\eta$ satisfying \cref{eq:defsmeta} with $\rhof=|\lambda|$, $\rhos=|\zeta|$ and $m=2r$, it holds \modr{$\exp(|R_{s}(p_m,q_r,\xi)|^2)< 1$}, i.e. the $\mSKROCK$ method is mean-square stable.
\end{theorem}
\begin{proof}
In \cite[Theorem 3.2]{AAV18} it is shown that if \modr{$p_m+\frac{1}{2}|q_r|^2<0$} and $|p_m|\leq \beta s^2$ then \modr{$\exp(|R_{s}(p_m,q_r,\xi)|^2)<1$}. 
We start noting that
\begin{equation}
\modr{p_m+\frac{1}{2}|q_r|^2< 0} \qquad \Longleftrightarrow  \qquad \modr{\Ps_m(\eta\lambda)(\lambda+\zeta)+\frac{1}{2}\Psi_r(\eta\lambda)^2|\mu|^2<  0}.
\end{equation}
From $\eta|\lambda|\leq \beta m^2$ and \cref{lemma:relBPhi} follows $\Psi_r(\eta\lambda)^2\leq \Ps_m(\eta\lambda)$, using $(\lambda,\zeta,\mu)\in\mathcal{S}^{mMS}$ yields
\begin{equation}\label{eq:stabineq}
\Ps_m(\eta\lambda)(\lambda+\zeta)+\frac{1}{2}\Psi_r(\eta\lambda)^2|\mu|^2\leq \modr{ \Ps_m(\eta\lambda)\left(\lambda+\zeta+\frac{1}{2}|\mu|^2\right)< 0}.
\end{equation}
Furthermore, \cref{thm:stab_mrkc} implies $|p_m|\leq \beta s^2$ (see \cite{AGR20}).
Thus, \modr{$\exp(|R_{s}(p_m,q_r,\xi)|^2)< 1$}.
\end{proof}
Even though \cref{thm:stab_mskrock} is stated for $\varepsilon=0$, numerical evidences show that it is valid for any damping $\varepsilon>0$; indeed, \cref{thm:stab_mrkc}, \cref{lemma:relBPhi} and \cite[Theorem 3.2]{AAV18} hold for $\varepsilon>0$.

We see from \eqref{eq:stabineq} that the stability of the $\mSKROCK$ scheme relies on the inequality $\Psi_r(z)^2\leq\Ps_m(z)$, where $\Psi_r(z)$ and $\Ps_m(z)$ are polynomials associated to the modified RKC scheme \cref{eq:defve} and the standard RKC scheme \cref{eq:defbue}, respectively. If instead of \cref{eq:defve} a standard RKC scheme is used then $\Ps_r(z/2)^2\leq \Ps_m(z)$ is needed for mean-square stability but this condition does not hold; hence, a modified RKC scheme is needed.
%If for the integration of \eqref{eq:defv} we naively used a standard RKC method, instead of the modified RKC method \eqref{eq:defve}, then we would obtain $q_r=\Ps_r(z/2)\mu\tau^{1/2}$ for some $r$ satisfying the stability conditions of \eqref{eq:defv} (instead of $q_r=\Psi_r(z)\mu\tau^{1/2}$ with $2r=m$). Hence, stability of the scheme would hold only if $\Ps_r(z/2)^2\leq \Ps_m(z)$ was true (as for $\varphi(z)$ in \cref{lemma:phisq}). However, such a relation is not satisfied and a modified RKC method must be used in order to replace $\Ps_r(z/2)$ by a different polynomial $\Psi_r(z)$ satisfying $\Psi_r(z)^2\leq\Ps_m(z)$. 

\subsection{Convergence analysis}\label{sec:accanalysis}
In this section we prove that the $\mSKROCK$ method has strong order $1/2$ and weak order $1$. We denote by $C_p^4(\Rn,\Rn)$ the space of functions from $\Rn$ to $\Rn$ four times continuously differentiable having derivatives with at most polynomial growth. \modr{For simplicity we suppose that $g$ is a vector valued function, the proofs for a matrix valued diffusion terms are similar.}
We start the convergence analysis stating a technical lemma, whose proof can be found in \cite[Section 4.5.2]{Ros20}.

%\begin{lemma}\label{lemma:bfbgerr}
%	There exists $C>0$ such that 
%	\begin{align}\label{eq:lipbfge}
%	|\bfe(x)-\bfe(y)|+|\bge(x)-\bge(y)|&\leq  C|x-y|, \\ \label{eq:errbfge}
%	|\bfe(x)-f(x)|+|\bge(x)-g(x)|&\leq  C(1+|x|)\eta
%	\end{align}
%	for all $x,y\in\Rn$.
%\end{lemma}
%
%\begin{lemma}\label{lemma:bfgereg}
%	If $\ff,\fs,g\in C_p^4(\Rb)$ then $\bfe,\bge\in C_p^4(\Rb)$.
%\end{lemma}
%
%\begin{lemma}\label{lemma:bnddiffKj}
%	The stages $K_j$ and $\overline Q_\eta$ of \cref{eq:mskrock} satisfy the estimate
%	\begin{align}\label{eq:bnddiffKj}
%	|\overline Q_\eta|+	|K_j-X_n|&\leq C(1+|X_n|)(\tau+|\Delta W_n|), \\ \label{eq:bndexpdiffKj}
%	|\exp(\overline Q_\eta|X_n)|+ |\exp(K_{j}-X_n|X_n)|&\leq C (1+|X_n|) \tau
%	\end{align}
%	for $j=1,\ldots,s$.
%\end{lemma}

\begin{lemma}\label{lemma:technical}
	Let $\ff,\fs,g$ be Lipschitz continuous, then there exists $C>0$ such that 
	\begin{align}\label{eq:lipbfge}
		|\bfe(x)-\bfe(y)|+|\bge(x)-\bge(y)|&\leq  C|x-y|, \\ \label{eq:errbfge}
		|\bfe(x)-f(x)|+|\bge(x)-g(x)|&\leq  C(1+|x|)\eta
	\end{align}
	for all $x,y\in\Rn$. Furthermore, the stages $K_j$ and $\overline Q_\eta$ of \cref{eq:mskrock} satisfy the estimate
		\begin{align}\label{eq:bnddiffKj}
		|\overline Q_\eta|+	|K_j-X_n|&\leq C(1+|X_n|)(\tau+|\Delta W_n|), \\ \label{eq:bndexpdiffKj}
		|\exp(\overline Q_\eta|X_n)|+ |\exp(K_{j}-X_n|X_n)|&\leq C (1+|X_n|) \tau
	\end{align}
	for $j=1,\ldots,s$.
\end{lemma}
\begin{lemma}\label{lemma:exprXnpu}
	Let $\ff,\fs,g$ be Lipschitz continuous, then the solution $X_{n+1}$ of \cref{eq:mskrock} satisfies
	\begin{equation}\label{eq:exprXnpu}
	X_{n+1} = X_n+\tau f(X_n)+g(X_n)\Delta W_n+R,
	\end{equation}
	with 
	\begin{equation}\label{eq:boundER}
	|\exp(R|X_n)|\leq C(1+|X_n|)\tau^{3/2} \qquad\mbox{and}\qquad \exp(|R|^2|X_n)^{1/2}\leq C(1+|X_n|)\tau^{3/2}.
	\end{equation}
	If, furthermore, $\ff,\fs\in C_p^1(\Rn,\Rn)$, then 
	$|\exp(R|X_n)|\leq C(1+|X_n|^q)\tau^{2}$, 
%	\begin{equation}\label{eq:imprboundEr}
%		|\exp(R|X_n)|\leq C(1+|X_n|^q)\tau^{2},
%	\end{equation}
	with $q\in\Nb$.
\end{lemma}
\begin{proof}
	It is shown recursively (see \cite[Lemma 4.19]{Ros20}) that
	\begin{equation}
	X_{n+1}	= X_n+ \frac{b_s}{b_1}U_{s-1}(\oz)\kappa_1\bQe +\tau \sum_{k=1}^s\frac{b_s}{b_k}U_{s-k}(\oz)\mu_k\bfe(\tilde K_{k-1}),
	\end{equation}
	where $\tilde K_0=K_0+\nu_1 \bQe$ and $\tilde K_k=K_k$ for $k=1,\ldots,s-1$. Since $\frac{b_s}{b_1}U_{s-1}(\oz)\kappa_1=1$, $\sum_{k=1}^s\frac{b_s}{b_k}U_{s-k}(\oz)\mu_k=1$ and $\bQe=\bge(X_n)\Delta W_n$, we can write
	\begin{equation}
	X_{n+1}= X_n+\tau f(X_n)+g(X_n)\Delta W_n+ R,
	\end{equation}
	with
	\begin{align}
	R&=R_1+R_2+R_3, & R_1&= \tau(\bfe(X_n)-f(X_n)),\\
	R_2&= (\bge(X_n)-g(X_n))\Delta W_n, &
	R_3&= \tau\sum_{k=1}^s\frac{b_s}{b_k}U_{s-k}(\oz)\mu_k(\bfe(\tilde K_{k-1})-\bfe(X_n)).
	\end{align}
	From  \cref{eq:errbfge},
	\begin{equation}\label{eq:R12sq}
	%\begin{split}
	|R_1|^2 \leq  C(1+|X_n|)^2\eta^2\tau^2, \qquad\qquad
	|R_2|^2 \leq  C(1+|X_n|)^2\eta^2\Delta W_n^2. 
	%\end{split}
	\end{equation}
	Since $\frac{b_s}{b_k}U_{s-k}(\oz)\mu_k\geq 0$ and $\sum_{k=1}^s\frac{b_s}{b_k}U_{s-k}(\oz)\mu_k=1$, using \cref{eq:lipbfge,eq:bnddiffKj} we obtain
	\begin{align}\label{eq:R3sq}
	\begin{split}
	|R_3|^2 &\leq   \tau^2  \max_{k=1,\ldots,s}|\bfe(\tilde K_{k-1})-\bfe(X_n)|^2\leq  C \tau^2 \max_{k=1,\ldots,s}|\tilde K_{k-1}-X_n|^2 \\
	&\leq  C(1+|X_n|)^2\tau^2(\tau+|\Delta W_n|)^2.
	\end{split}
	\end{align}
	From $\exp(R_2|X_n)=0$, \cref{eq:R12sq,eq:R3sq}, using Jensen's inequality we get
	\begin{align}
	|\exp(R|X_n)|&\leq  C(1+|X_n|)(\eta+\tau^{1/2})\tau ,&
	\exp(|R|^2|X_n)&\leq C(1+|X_n|)^2(\eta^2\tau+\eta^2+\tau^{2})\tau.
	\end{align}
	This proves \cref{eq:boundER} using \cref{eq:defsmeta}, from where we deduce $\eta\leq 8\tau$ (as $\beta s^2\geq 1$ and $m\geq 2$).
	
	To show the improved estimate on $|\exp(R|X_n)|$ we suppose $\ff,\fs\in C_p^1(\Rb)$. It can be shown recursively that $\bfe\in C_p^1(\Rb)$ and from \cref{eq:bndexpdiffKj} it holds $|\exp(\tilde K_{k-1}-X_n|X_n)|\leq C (1+|X_n|)\tau$, thus
	\begin{align}\label{eq:imprestR3}
	\begin{split}
	|\exp(R_3|X_n)|&\leq  \tau\sum_{k=1}^s\frac{b_s}{b_k}U_{s-k}(\oz)\mu_k|\exp(\bfe(\tilde K_{k-1})-\bfe(X_n)|X_n)|\\
	&\leq \tau\max_{k=1,\ldots,s}|\exp(\bfe(\tilde K_{k-1})-\bfe(X_n)|X_n)|\\
	&\leq C\tau (1+|X_n|^q) \max_{k=1,\ldots,s}|\exp(\tilde K_{k-1}-X_n|X_n)|
	\leq  C(1+|X_n|^{\modr{q+1}})\tau^2,
	\end{split}
	\end{align}
	where we used $\bfe\in C_p^1(\Rb)$ and \cref{eq:bnddiffKj} to bound the derivative of $\bfe$ in $[X_n,\tilde K_{k-1}]$ by $C(1+|X_n|^q)$.
 	Using $|\exp(R_1|X_n)|\leq C(1+|X_n|)\tau^2 $ and \cref{eq:imprestR3} yields $|\exp(R|X_n)|\leq C (1+|X_n|^{\modr{q+1}})\tau^2$. 
\end{proof}

\begin{theorem}\label{thm:conv}
	Consider the system of SDEs \cref{eq:msde} on $[0,T]$, $T>0$. Assume that $\ff,\fs,g\in C_p^4(\Rb)$ are Lipschitz continuous, then the $\mSKROCK$ method has strong order $1/2$ and weak order $1$, i.e.
	\begin{align}\label{eq:strongskrock}
	\exp(|X(\tn)- X_n|^2)^{1/2} &\leq C(1+\exp(|X_0|^2))^{1/2}\tau^{1/2}, \\ \label{eq:weakskrock}
	|\exp(\psi(X(\tn))-\exp(\psi(X_n))| &\leq C(1+\exp(|X_0|^q))\tau  
	\end{align}
	for $t_n=n\tau \leq T$ and all $\psi\in C_p^4(\Rb)$, where $C$ is independent from $n,\tau$.
\end{theorem}
\begin{proof}
	As $\ff,\fs,g$ are Lipschitz continuous, \modr{doing a stochastic Taylor expansion of} \cref{eq:msde} with initial value $X(\tn)=X_n$, we obtain
	\begin{equation}
	X(\tnpu)=X_n+\tau f(X_n)+g(X_n)\Delta W_n+\overline{R}
	\end{equation}
	with $|\exp(\overline R|X_n)|\leq C(1+|X_n|)\tau^{\modr{2}}$ and $\exp(|\overline R|^2|X_n)^{1/2}\leq C(1+|X_n|)\tau$. Therefore, it follows from \cref{lemma:exprXnpu} that the local errors satisfy
	\begin{align}
	|\exp(X_{n+1}-X(\tnpu)|X_n)|&\leq C(1+|X_n|)\tau^{3/2}, &
	\exp(|X_{n+1}-X(\tnpu)|^2|X_n)^{1/2}&\leq C(1+|X_n|)\tau.
	\end{align}
	The classical result \cite[Theorem 1.1]{MiT03}, which asserts the global order of convergence from the local error, implies estimate \cref{eq:strongskrock}.
	From \cref{lemma:exprXnpu} and the Itô formula we obtain the local error estimate
	\begin{equation}
	|\exp(\psi(X(\tnpu))-\psi(X_{n+1})|X_n)|\leq  C(1+|X_n|^q)\tau^2.
	\end{equation}
	Next we need to show that the moments $\exp(|X_n|^{2k})$ are bounded for $k\in\Nb$ and all $n,\tau$ with $0\leq n\tau\leq T$ uniformly with respect to all small enough $\tau$. Using \cite[Lemma 2.2]{MiT03} this follows from \cref{eq:bnddiffKj},\cref{eq:bndexpdiffKj}. Finally from the local error estimate, the bounded moments and the regularity assumption on $\ff,\fs$ and $g$ we obtain \cref{eq:weakskrock} from the classical result for weak convergence \cite[Theorem 2.1]{MiT03}.
\end{proof}

\section{Numerical experiments}\label{sec:numexp}
Through a series of numerical experiments, we illustrate here the accuracy of the $\mSKROCK$ method of \cref{sec:mskrockalg} and compare its computational cost against the cost of the standard SK-ROCK scheme; which is given by \cref{eq:mskrock} but where $\bfe,\bge$ are replaced by $f,g$ and the stability condition is $\tau\rho\leq \beta s^2$, with $\rho$ the spectral radius of the Jacobian of $f$. At first, we confirm the strong and weak convergence properties of the $\mSKROCK$ scheme on a nonstiff problem, where we fix the number of stages beforehand. Then we do the same but on a stiff problem, letting the scheme automatically choose the number of stages based on the spectral radii and the step size. \modr{To do so, at each time step the spectral radii of the Jacobians of $\ff,\fs$ are estimated employing a cheap nonlinear power method \cite{Lin72,Ver80}, then the number of stages $s,m$ are chosen according to \cref{eq:defsmeta}.}
For the last two examples we consider the application of the $\mSKROCK$ method to more challenging problems, first on a chemical Langevin equation and then on a stochastic heat equation with multiplicative noise. The last experiment has been performed with the help of the C++ library \texttt{libMesh} \cite{KPS06}.

We note that while we compare the $\mSKROCK$ method only to the SK-ROCK method, reference \cite{AAV18} contains comparisons of SK-ROCK with many other stabilized methods (S-ROCK, S-ROCK2, PSK-ROCK) and for the type of problems considered here SK-ROCK shows the best performance.

\subsection{Nonstiff problem convergence experiment}\label{sec:nsconv}
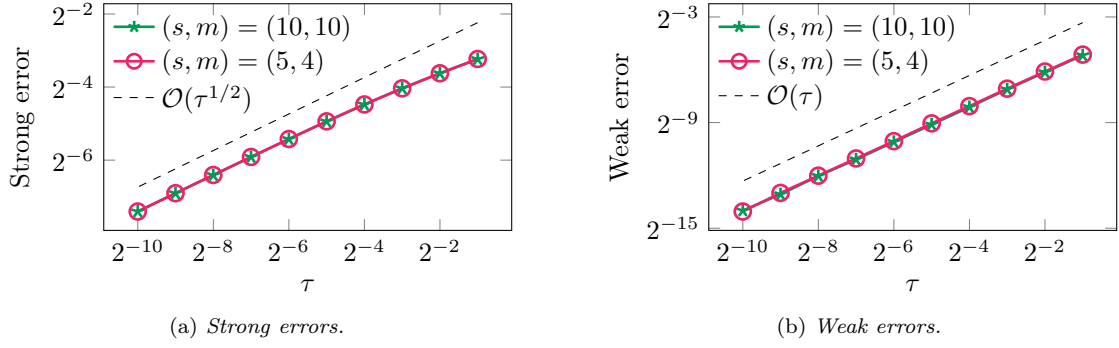
\begin{figure}
	\begin{center}
		\begin{subfigure}[t]{\subfigsize\textwidth}
			\centering
			\begin{tikzpicture}[scale=\plotimscale]
			\begin{loglogaxis}[height=\aspectratio*\plotimsized\textwidth,width=\plotimsized\textwidth,legend columns=1,legend style={draw=\legendboxdraw,fill=\legendboxfill,at={(0,1)},anchor=north west},log basis x={2},log basis y={2},legend cell align={left},
			xlabel={$\tau$}, ylabel={Strong error},label style={font=\normalsize},tick label style={font=\normalsize},legend image post style={scale=\legendmarkscale},legend style={nodes={scale=\legendfontscale, transform shape}},]
			\addplot[color=colorone,solid,line width=\plotlinewidth pt,mark=\markone,mark size=\plotmarksizeu pt] table [x=dt,y=errstrong,col sep=comma] 
			{data/exp_conv/SK-mROCK_s_10_m_10.csv};\addlegendentry{$(s,m)=(10,10)$}
			\addplot[color=colorthree,line width=\plotlinewidth pt,mark=\markthree,mark size=\plotmarksized pt] table [x=dt,y=errstrong,col sep=comma] 
			{data/exp_conv/SK-mROCK_s_5_m_4.csv};\addlegendentry{$(s,m)=(5,4)$}
			\addplot[black,dashed,domain=1e-3:0.5] (x,0.3*x^0.5);\addlegendentry{$\bigo{\tau^{1/2}}$}
			\end{loglogaxis}
			\end{tikzpicture}
			\caption{Strong errors.}
%			\label{fig:strongerr}
		\end{subfigure}\hfill
		\begin{subfigure}[t]{\subfigsize\textwidth}
			\centering
			\begin{tikzpicture}[scale=\plotimscale]
			\begin{loglogaxis}[height=\aspectratio*\plotimsized\textwidth,width=\plotimsized\textwidth,legend columns=1,legend style={draw=\legendboxdraw,fill=\legendboxfill,at={(0,1)},anchor=north west},log basis x={2},log basis y={2},legend cell align={left},
			xlabel={$\tau$}, ylabel={Weak error},label style={font=\normalsize},tick label style={font=\normalsize},legend image post style={scale=\legendmarkscale},legend style={nodes={scale=\legendfontscale, transform shape}},]
			\addplot[color=colorone,solid,line width=\plotlinewidth pt,mark=\markone,mark size=\plotmarksizeu pt] table [x=dt,y=errweak,col sep=comma] 
			{data/exp_conv/SK-mROCK_s_10_m_10.csv};\addlegendentry{$(s,m)=(10,10)$}
			\addplot[color=colorthree,line width=\plotlinewidth pt,mark=\markthree,mark size=\plotmarksized pt] table [x=dt,y=errweak,col sep=comma] 
			{data/exp_conv/SK-mROCK_s_5_m_4.csv};\addlegendentry{$(s,m)=(5,4)$}
			\addplot[black,dashed,domain=1e-3:0.5] (x,0.2*x);\addlegendentry{$\bigo{\tau}$}
			\end{loglogaxis}
			\end{tikzpicture}
			\caption{Weak errors.}
%			\label{fig:weakerr}
		\end{subfigure}
	\end{center}
	\caption{Nonstiff convergence experiment. Strong and weak errors of $\mSKROCK$ vs. step size $\tau$, for different stage choices.}
	\label{fig:convmskrock}
\end{figure}

We perform a convergence experiment on the following SDE, taken from \cite{AAV18},
\begin{equation}
	\dif X(t) = \left(\frac{1}{4}X(t)+\frac{1}{2}\sqrt{X(t)^2+1}\right)\dif t+\sqrt{\frac{X(t)^2+1}{2}}\dif W(t), \qquad\qquad X(0)=0,
\end{equation}
where the exact solution is $X(t)=\sinh(\frac{t}{2}+\frac{W(t)}{\sqrt{2}})$. We let $\ff(X)=\frac{1}{2}\sqrt{X^2+1}$ and $\fs(X)=\frac{1}{4}X$. Considering the step sizes $\tau=2^{-k}$, for $k=1,\ldots,10$, we display the strong $\exp(|X(T)-X_N|^2)^{1/2}$ and weak $|\exp(\asinh(X(T)))-\exp(\asinh(X_N))|$ errors at time $T=1=N\tau$ in \cref{fig:convmskrock}, using $10^6$ samples and $(s,m)=(5,4)$ or $(s,m)=(10,10)$. We observe that the method converges with the predicted orders of accuracy and the error is essentially independent of the stages number.

\subsection{Multiscale problem convergence experiment}\label{sec:sconv}
We consider a chemical Langevin model of dimerization reactions in a genetic network \cite{BHJ03}. The model consists of 7 species and 10 reactions, described by the equations
\begin{equation}\label{eq:cle}
	\dif X(t)= \sum_{j=1}^{l} \nu_j a_j(X(t))\dif t+ \sum_{j=1}^{l} \nu_j\sqrt{a_j(X(t))}\dif W_j(t)\quad t\in(0,T], \qquad\qquad X(0)= X_0,
\end{equation}
where $l=10$, $X(t)\in\Rb^{7}$, $T=10$ and $\nu_j$, $a_j(x)$ are derived from the chemical reaction system introduced in \cite{BHJ03}.  We consider the same initial conditions as in  \cite{BHJ03} but multiplied by $10^3$.

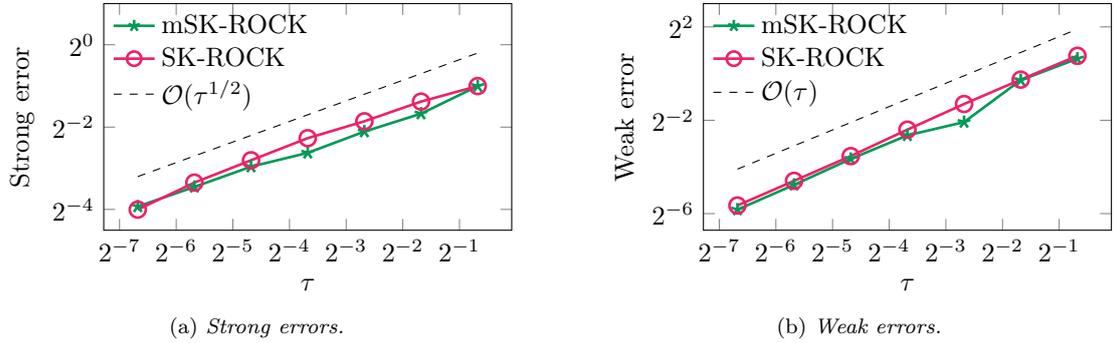
\begin{figure}
	\begin{center}
		\begin{subfigure}[t]{\subfigsize\textwidth}
			\centering
			\begin{tikzpicture}[scale=\plotimscale]
			\begin{loglogaxis}[height=\aspectratio*\plotimsized\textwidth,width=\plotimsized\textwidth,legend columns=1,ymax=2,legend style={draw=\legendboxdraw,fill=\legendboxfill,at={(0,1)},anchor=north west},log basis x={2},log basis y={2},legend cell align={left},
			xlabel={$\tau$}, ylabel={Strong error},label style={font=\normalsize},tick label style={font=\normalsize},legend image post style={scale=\legendmarkscale},legend style={nodes={scale=\legendfontscale, transform shape}},]
			\addplot[color=localcolor,solid,line width=\plotlinewidth pt,mark=\localmark,mark size=\plotmarksizeu pt,select coords between index={0}{6}] table [x=dt,y=errstrong,col sep=comma] 
			{data/exp_conv_stiff/SK-mROCK.csv};\addlegendentry{$\mSKROCK$}
			\addplot[color=classicalcolor,line width=\plotlinewidth pt,mark=\classicalmark,mark size=\plotmarksized pt,select coords between index={0}{6}] table [x=dt,y=errstrong,col sep=comma] 
			{data/exp_conv_stiff/SK-ROCK.csv};\addlegendentry{SK-ROCK}
			\addplot[black,dashed,domain=0.009765625:0.625] (x,1.1*x^0.5);\addlegendentry{$\bigo{\tau^{1/2}}$}
			\end{loglogaxis}
			\end{tikzpicture}
			\caption{Strong errors.}
%			\label{fig:strongerrstiff}
		\end{subfigure}\hfill
		\begin{subfigure}[t]{\subfigsize\textwidth}
			\centering
			\begin{tikzpicture}[scale=\plotimscale]
			\begin{loglogaxis}[height=\aspectratio*\plotimsized\textwidth,width=\plotimsized\textwidth,legend columns=1,ymax=8,legend style={draw=\legendboxdraw,fill=\legendboxfill,at={(0,1)},anchor=north west},log basis x={2},log basis y={2},legend cell align={left},
			xlabel={$\tau$}, ylabel={Weak error},label style={font=\normalsize},tick label style={font=\normalsize},legend image post style={scale=\legendmarkscale},legend style={nodes={scale=\legendfontscale, transform shape}},]
			\addplot[color=localcolor,line width=\plotlinewidth pt,mark=\localmark,mark size=\plotmarksizeu pt,select coords between index={0}{6}] table [x=dt,y=errweak,col sep=comma] 
			{data/exp_conv_stiff/SK-mROCK.csv};\addlegendentry{$\mSKROCK$}
			\addplot[color=classicalcolor,line width=\plotlinewidth pt,mark=\classicalmark,mark size=\plotmarksized pt,select coords between index={0}{6}] table [x=dt,y=errweak,col sep=comma] 
			{data/exp_conv_stiff/SK-ROCK.csv};\addlegendentry{SK-ROCK}
			\addplot[black,dashed,domain=0.009765625:0.625] (x,6*x);\addlegendentry{$\bigo{\tau}$}
			\end{loglogaxis}
			\end{tikzpicture}
			\caption{Weak errors.}
%			\label{fig:weakerrstiff}
		\end{subfigure}
	\end{center}
	\caption{Multiscale convergence experiment. Strong and weak errors of $\mSKROCK$ and SK-ROCK vs. step size $\tau$.}
	\label{fig:sconv}
\end{figure}

We order the reaction terms $\nu_j a_j(x)$ from the fastest to the slowest (the sequence $\rho_j$ of the spectral radii of the Jacobians of $\nu_j a_j(x)$, evaluated on a typical path $X(t)$, is decreasing) and let
\begin{equation}
	\ff(x)=\sum_{j=1}^{3} \nu_j a_j(x), \qquad\qquad \fs(x)=\sum_{j=4}^{10} \nu_j a_j(x),
\end{equation}
hence $\ff$ represents the three fastest reactions.
We run the $\mSKROCK$ method over $10^5$ Brownian paths with step size $\tau=T/2^k$ for $k=4,\ldots,10$ and measure the strong and weak errors committed against reference solutions computed on the same paths but using the SK-ROCK method with a step size $\tau=T/2^{12}$. As weak error we consider the error committed on the second moment of $X_6$. Differently from \cref{sec:nsconv} we let the $\mSKROCK$ method automatically choose the number of stages $s,m$. We observe in \cref{fig:sconv} that the $\mSKROCK$ method converges with the right orders and have similar errors as the SK-ROCK scheme.

\subsection{E. Coli bacteria heat shock response}
We consider a chemical Langevin equation modeling E. coli bacteria's protein denaturation under heat shocks. The original deterministic model is introduced in \cite{KSY01}, while in \cite{CLP04,HAL12} it is considered as a chemical reaction system. 

The model consists of 28 species and 61 reactions, it is described by \cref{eq:cle} with $l=61$ and $X(t)\in\Rb^{28}$. The initial condition is the same as in \cite{HAL12} but multiplied by 100 and we let $T=10$. The parameters $\nu_j$, $a_j(x)$ are derived from the chemical reactions described in \cite[Section 7.2]{HAL12} and the terms $\nu_j a_j(x)$ are ordered from the fastest to the slowest as explained in \cref{sec:sconv}. For $r=0,\ldots,10$ we define 
\begin{equation}
	\ff^r(x)=\sum_{j=1}^{r} \nu_j a_j(x), \qquad\qquad \fs^r(x)=\sum_{j=r+1}^{61} \nu_j a_j(x),
\end{equation}
hence $\ff^r$ is defined by the $r$ fastest reactions and $\fs^r$ by the remaining ones. Observe that for $r=0$ it holds $\ff^0=0$ and thus all the reactions are considered to be slow.

\begin{figure}
	\begin{center}
		\begin{subfigure}[t]{\subfigsize\textwidth}
			\centering
			\begin{tikzpicture}[scale=\plotimscale]
			\begin{semilogyaxis}[height=\aspectratio*\plotimsized\textwidth,width=\plotimsized\textwidth,legend columns=1,legend style={draw=\legendboxdraw,fill=\legendboxfill,at={(0,0)},anchor=south west},log basis y={10},ymax=40000000,legend cell align={left},
			xlabel={$r$}, ylabel={Spectral radii},label style={font=\normalsize},tick label style={font=\normalsize},legend image post style={scale=\legendmarkscale},legend style={nodes={scale=\legendfontscale, transform shape}},]
			\addplot[color=localcolorF,solid,line width=\plotlinewidth pt,mark=\localmarkF,mark size=\plotmarksizeu pt,select coords between index={1}{10}] table [x=Nreac,y=rF,col sep=comma] 
			{data/ecoli/data.csv};\addlegendentry{$\rhof$}
			\addplot[color=localcolorS,line width=\plotlinewidth pt,mark=\localmarkS,mark size=\plotmarksizeu pt] table [x=Nreac,y=rS,col sep=comma] 
			{data/ecoli/data.csv};\addlegendentry{$\rhos$}
			\end{semilogyaxis}
			\end{tikzpicture}
			\caption{Spectral radii of the Jacobians of $\fs^r$ and $\ff^r$.}
			\label{fig:ecolirho}
		\end{subfigure}\hfill
		\begin{subfigure}[t]{\subfigsize\textwidth}
			\centering
			\begin{tikzpicture}[scale=\plotimscale]
			\begin{semilogyaxis}[height=\aspectratio*\plotimsized\textwidth,width=\plotimsized\textwidth,legend columns=1,legend style={draw=\legendboxdraw,fill=\legendboxfill,at={(0,0)},anchor=south west},log basis y={10},legend cell align={left},ymin=1,
			xlabel={$r$}, ylabel={Stages},label style={font=\normalsize},tick label style={font=\normalsize},legend image post style={scale=\legendmarkscale},legend style={nodes={scale=\legendfontscale, transform shape}},]
			\addplot[color=localcolorF,line width=\plotlinewidth pt,mark=\localmarkF,mark size=\plotmarksizeu pt,select coords between index={1}{10}] table [x=Nreac,y=m,col sep=comma] 
			{data/ecoli/data.csv};\addlegendentry{$m$}
			\addplot[color=localcolorS,solid,line width=\plotlinewidth pt,mark=\localmarkS,mark size=\plotmarksizeu pt] table [x=Nreac,y=s,col sep=comma] 
			{data/ecoli/data.csv};\addlegendentry{$s$}
			\end{semilogyaxis}
			\end{tikzpicture}
			\caption{Number of stages $s$, $m$.}
			\label{fig:ecolism}
		\end{subfigure}
	\end{center}
	\caption{E. Coli experiment. Spectral radii and number of stages vs. the number $r$ of fast reactions put into $\ff^r$.}
	\label{fig:ecolirhosm}
\end{figure}
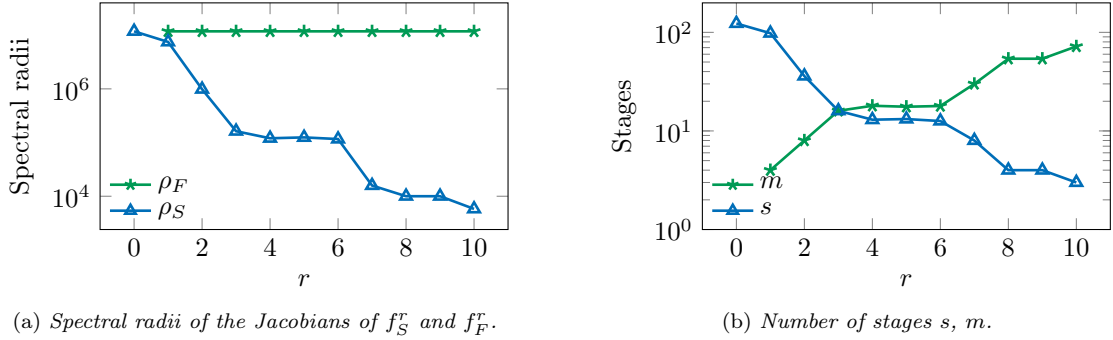

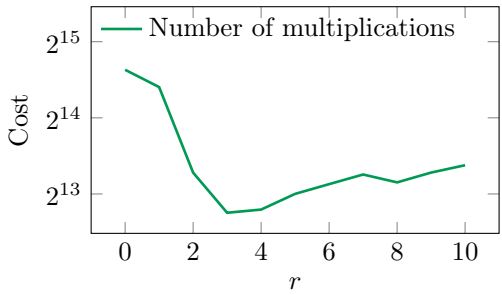
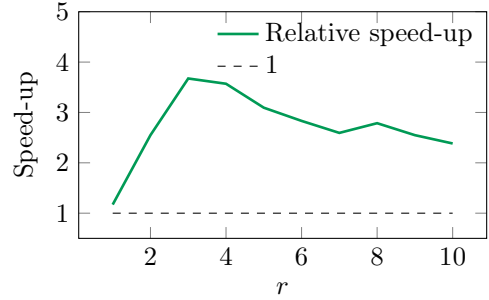
\begin{figure}
	\begin{center}
		\begin{subfigure}[t]{\subfigsize\textwidth}
			\centering
			\begin{tikzpicture}[scale=\plotimscale]
			\begin{semilogyaxis}[height=\aspectratio*\plotimsized\textwidth,width=\plotimsized\textwidth,legend columns=1,legend style={draw=\legendboxdraw,fill=\legendboxfill,at={(0,1)},anchor=north west},log basis y={2},ymax=45000,legend cell align={left},
			xlabel={$r$}, ylabel={Cost},label style={font=\normalsize},tick label style={font=\normalsize},legend image post style={scale=\legendmarkscale},legend style={nodes={scale=\legendfontscale, transform shape}},]
			\addplot[color=colorone,solid,line width=\plotlinewidth pt,mark=none,mark size=\plotmarksizeu pt] table [x=Nreac,y=cost,col sep=comma] 
			{data/ecoli/data.csv};\addlegendentry{Number of multiplications}
			\end{semilogyaxis}
			\end{tikzpicture}
			\caption{Cost in number of multiplications.}
			\label{fig:ecolicost}
		\end{subfigure}\hfill
		\begin{subfigure}[t]{\subfigsize\textwidth}
			\centering
			\begin{tikzpicture}[scale=\plotimscale]
			\begin{axis}[height=\aspectratio*\plotimsized\textwidth,width=\plotimsized\textwidth,legend columns=1,legend style={draw=\legendboxdraw,fill=\legendboxfill,at={(1,1)},anchor=north east},ymin=0.5,ymax=5,legend cell align={left},
			xlabel={$r$}, ylabel={Speed-up},label style={font=\normalsize},tick label style={font=\normalsize},legend image post style={scale=\legendmarkscale},legend style={nodes={scale=\legendfontscale, transform shape}},ytick={1,2,3,4,5}]
			\addplot[color=colorone,solid,line width=\plotlinewidth pt,mark=none,mark size=\plotmarksizeu pt,select coords between index={1}{10}] table [x=Nreac,y=speed,col sep=comma] 
			{data/ecoli/data.csv};\addlegendentry{Relative speed-up}
			\addplot[black,dashed,domain=1:10] (x,1);\addlegendentry{$1$}
			\end{axis}
			\end{tikzpicture}
			\caption{Relative speed-up with respect to SK-ROCK.}
			\label{fig:ecolispeedup}
		\end{subfigure}
	\end{center}
	\caption{E. Coli experiment. Cost and relative speed-up of $\mSKROCK$.}
	\label{fig:efficiency}
\end{figure}

Let $\tau=T/2^{12}$ be fixed, for each value of $r=0,\ldots,10$ we run the $\mSKROCK$ scheme and measure the following data: the mean values of $\rhof$, $\rhos$, $s$, $m$ along the integration interval and the code efficiency in terms of total multiplications needed to evaluate $\ff^r$ and $\fs^r$. For $r=0$ we have $\ff^0=0$ and thus the original SK-ROCK scheme is used with $f=\fs^0$. We display in \cref{fig:ecolirho,fig:ecolism} the values of $\rhof$, $\rhos$ and $s$, $m$, respectively. We see how $\rhos$ decreases as $r$ increases, indeed more fast reactions are put into $\ff^r$, as a consequence $s$ decreases as well. In order to compensate the decreasing stabilization made by the ``outer'' scheme, the ``inner'' method must increase the number of stages $m$, see \cref{fig:ecolism}.
In \cref{fig:ecolicost} we show the cost of the scheme, defined as the total number of multiplications needed by $\mSKROCK$ in order to evaluate $\bfe$ and $\bge$. For $r=0$ we have the cost of SK-ROCK and for $r=1,\ldots,10$ the cost of $\mSKROCK$. In \cref{fig:ecolispeedup} we show the relative speed-up of $\mSKROCK$ with respect to SK-ROCK, defined as the cost of SK-ROCK ($r=0$ in \cref{fig:ecolicost}) divided by the cost of $\mSKROCK$ for $r=1,\ldots,10$. We note that the speed-up reaches a maximal value and then decreases as more terms are put into $\ff^r$ and thus its evaluation becomes more expensive.

\subsection{Diffusion across a narrow channel with multiplicative space-time noise}\label{sec:channel}
Here, we consider a stochastic heat equation with multiplicative noise defined on a domain which requires local mesh refinement. We compare the efficiency of the $\mSKROCK$ and SK-ROCK method as the geometry imposes increasingly severe stability constraints. This problem is a stochastic version of a PDE problem studied in \cite{AGR20}.

\begin{figure}
	\begin{center}
		\begin{subfigure}[t]{\subfigsizemed\textwidth}
			\centering
			\begin{tikzpicture}[spy using outlines= {circle, magnification=6, connect spies}]
			\coordinate (a) at (-2.19,0);
			\coordinate (b) at (1,0);
			\coordinate (spypoint) at (b);
			\coordinate (magnifyglass) at (2.2,-1.2);
			\node at (a) {\includegraphics[width=0.16\textwidth]{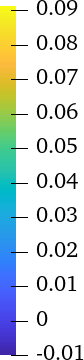}};
			\node at (b) {\includegraphics[width=0.74\textwidth]{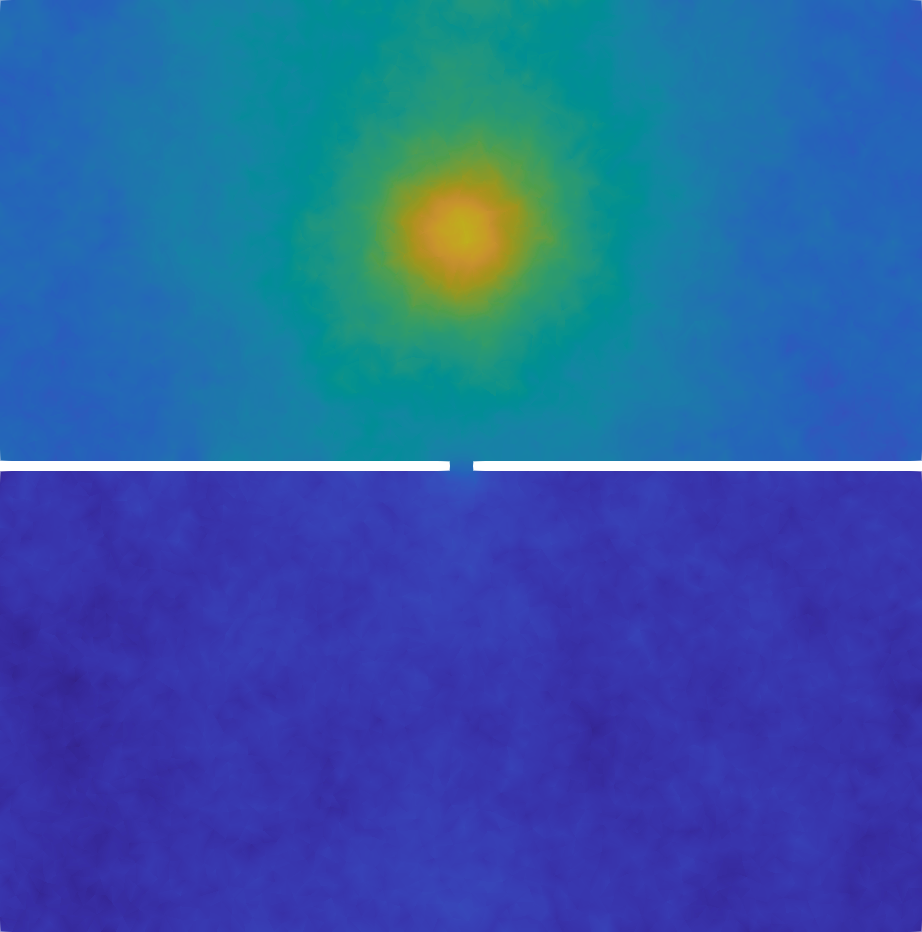}};
			\spy [BrickRed, size=2cm,magnification=6] on (spypoint) in node[fill=white] at (magnifyglass);
			\end{tikzpicture}
			\caption{Solution for $\delta=1/2^2$.}
			\label{fig:eff2_sto_sola}
		\end{subfigure}\hfill
		\begin{subfigure}[t]{\subfigsizemed\textwidth}
			\centering
			\begin{tikzpicture}[spy using outlines= {circle, magnification=6, connect spies}]
			\coordinate (a) at (-2.19,0);
			\coordinate (b) at (1,0);
			\coordinate (spypoint) at (b);
			\coordinate (magnifyglass) at (2.2,-1.2);
			\node at (a) {\includegraphics[width=0.16\textwidth]{images/exp_eff2/bar.png}};
			\node at (b) {\includegraphics[width=0.74\textwidth]{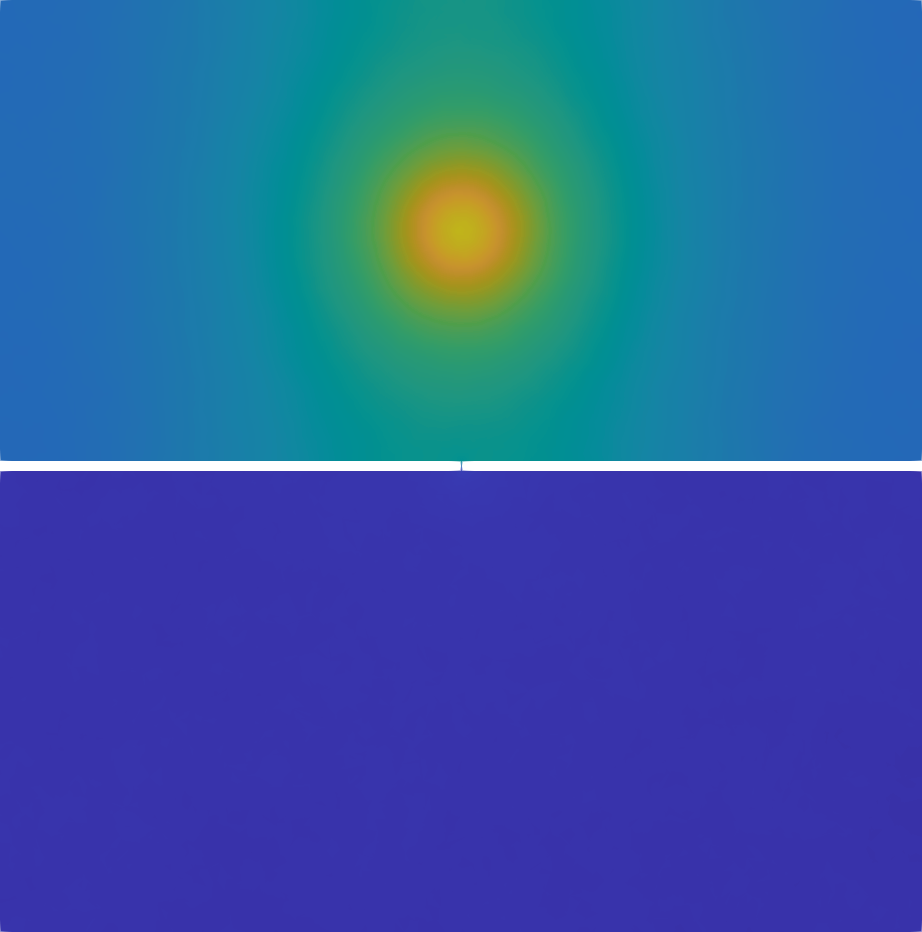}};
			\spy [BrickRed, size=2cm,magnification=6] on (spypoint) in node[fill=white] at (magnifyglass);
			\end{tikzpicture}
			\caption{Solution for $\delta=1/2^7$.}
			\label{fig:eff2_sto_solb}
		\end{subfigure}
	\end{center}
	\caption{Narrow channel. Numerical solutions at $t=10$ for a channel width $\delta=1/2^2$ (left) or $\delta=1/2^7$ (right).}
	\label{fig:eff2_sto_sol}
\end{figure}

We consider the next heat equation with multiplicative noise, colored in space and white in time:
\begin{equation}\label{eq:stonarch}
	\begin{aligned}
		\dif u &=(\Delta u+b)\dif t+G(u)\dif W \qquad \text{in }\Omega_\delta\times [0,T],\\ 
		\nabla u\cdot \bm{n}&=0 \hspace{116pt} \text{in } \partial\Omega_\delta \times [0,T],\\ 
		u&=0 \hspace{116pt} \text{in } \Omega_\delta\times\{0\},
	\end{aligned}
\end{equation}
where $T=0.1$ and $\Omega_\delta$ is a domain consisting in two $10\times 5$ rectangles linked together by a narrow channel $\delta\times 0.05$ of width $\delta>0$, see \cref{fig:eff2_sto_sol}. The source term $b(\bx,t)=\sin(10\pi t)^2e^{-5\nld{\bx-\bf{c}}^2}$ is a Gaussian centered in $\bf{c}$, the center of the upper rectangle in $\Omega_\delta$. We define $G:L^2(\Omega_\delta)\rightarrow \mathcal{L}(L^2(\Omega_\delta),L^2(\Omega_\delta))$ by $G(u)(v)(\bx)=u(\bx)v(\bx)$ and $W(t)$ is a $Q$-Wiener process defined by a covariance operator $Q:L^2(\Omega_\delta)\rightarrow L^2(\Omega_\delta)$, i.e. $W(t)$ satisfies 
\begin{equation}\label{eq:WQ}
	\exp(\langle W(t),h\rangle)=0 \qquad \mbox{and} \qquad \exp(\langle W(t),h\rangle^2)= t\langle Qh,h\rangle
\end{equation}
for all $h\in L^2(\Omega_\delta)$, where $\langle\cdot,\cdot,\rangle$ is the inner product in $L^2(\Omega_\delta)$. For $h\in L^2(\Omega_\delta)$ we define $Q$ by
\begin{equation}
	(Qh)(\bx)=\int_{\Omega_\delta}q(\bx-\by)h(\by)\dif\by, \qquad\mbox{where}\qquad q(\bx)=\frac{\alpha}{\pi}e^{-\alpha\nld{\bx}^2}
\end{equation}
is an approximation of the Dirac delta function and $\alpha=100$.

In $\Omega_\delta$, we define a Delaunay triangulation $\mathcal{M}$ composed by simplicial elements having maximal size $H\approx 0.015$. 
Let $V=\mbox{span}\{\varphi_i\,:\, i=1,\ldots,N\}$ be a first-order discontinuous Galerkin finite element (DG-FE) \cite{PiE12} space on $\mathcal{M}$ and $\Delta_H:V\rightarrow V$ the DG-FE discretization of the Laplacian. Then, the semidiscrete problem corresponding to \cref{eq:stonarch} is to find the process $u_H(t)=\sum_{i=1}^Nu_i(t)\varphi_i$ satisfying
\begin{equation}\label{eq:parsd}
	\dif u_H= (\Delta_H u_H+P_H b)\dif t+P_HG(u_H)\dif\widehat{W},
\end{equation}
where $P_H:L^2(\Omega_\delta)\rightarrow V$ is the orthogonal projection operator and $\widehat{W}(t)\in V$ is the numerical counterpart of $W(t)$ in \cref{eq:WQ}, hence it satisfies
\begin{equation}
	\exp(\langle \widehat{W}(t),h\rangle)=0 \qquad \mbox{and} \qquad \exp(\langle \widehat{W}(t),h\rangle^2)=t\langle Q h,h\rangle
\end{equation}
for all $h\in V$. We set
\begin{equation}
	\widehat{W}(t)=\sum_{i=1}^N\gamma_i^{1/2} e_i\beta_i(t),
\end{equation}
where $\{e_i\}_{i=1}^N$ is an orthonormal basis of $V$, $\gamma_i\geq 0$ for $i=1,\ldots,N$ and $\{\beta_i(t)\}_{i=1}^N$ is a sequence of independently and identically distributed Brownian motions. We have $\exp(\langle \widehat{W}(t),h\rangle)=0$ and since $\exp(\langle \widehat{W}(t),e_i\rangle^2)=t\gamma_i$ we set
\begin{equation}
	\gamma_i=\langle Qe_i,e_i\rangle=\int_{\Omega_\delta\times\Omega_\delta}q(\bx-\by)e_i(\bx)e_i(\by)\dif\bx\dif\by.
\end{equation}
Note that $\widehat{W}(t)$ is a $Q$-Wiener process in $V$ with covariance operator $Q_H$ defined by $Q_He_i=\gamma_i e_i$.

Taking the inner product on both sides of \cref{eq:parsd} with respect to $\varphi_j$ we obtain the equivalent equation
\begin{equation}\label{eq:parode}
	\dif X(t) = (AX(t)+M^{-1}\widehat b(t))\dif t+M^{-1}\widehat{G}(X(t))\dif B(t),
\end{equation}
with $X(t)=(u_i(t))_{i=1}^N$, $M$ the mass matrix, $A$ the stiffness matrix, $B(t)=(\beta_i(t))_{i=1}^N$ an $N$-dimensional Wiener process and $\widehat b(t)\in\Rb^N$, $\widehat G(X)\in \Rb^{N\times N}$ are defined by
\begin{equation}
	\widehat b_j(t)=\langle b(t),\varphi_j\rangle, \qquad\qquad \widehat G(X)_{ji}=\gamma_i^{1/2}\langle G(u_H)(e_i),\varphi_j\rangle.
\end{equation}
By setting
\begin{equation}
	f(t,X)=A X+M^{-1}\widehat b(t), \qquad\qquad g(X)=M^{-1}\widehat{G}(X)
\end{equation}
we obtain \cref{eq:msde}, in nonautonomous form. Note that the orthonormal basis $\{e_i\}_{i=1}^N$ can be computed locally on each element and $M$ is easy to invert since it is block-diagonal. Therefore, application of SK-ROCK to \cref{eq:parode} leads to a truly explicit method. 

\begin{figure}
	\begin{center}
		\begin{subfigure}[t]{\subfigsizemed\textwidth}
			\centering
			\includegraphics[trim=0cm 0cm 0cm 0cm, clip, width=0.24\textheight]{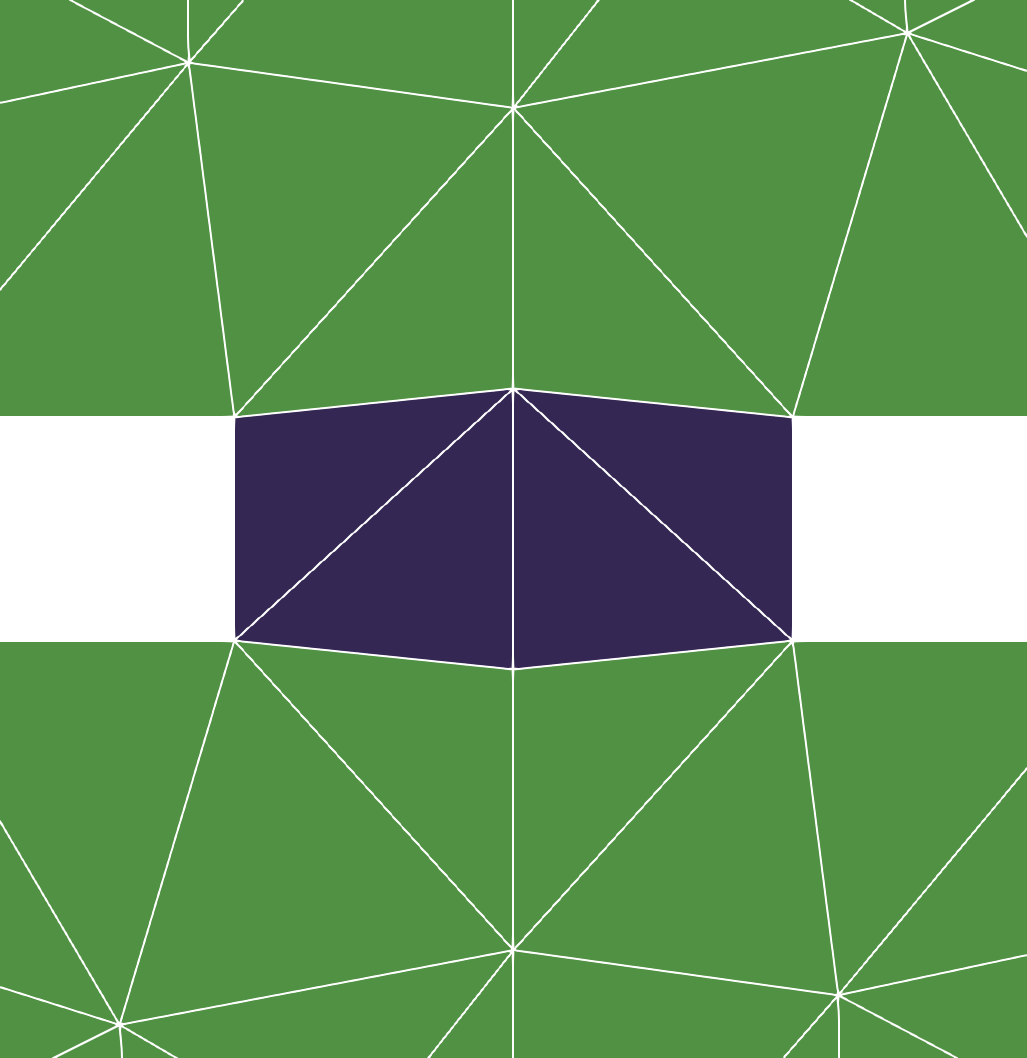}
			\caption{Zoom over the channel with $\delta=1/2^2$.}
			\label{fig:eff2_sto_grida}
		\end{subfigure}\hfill
		\begin{subfigure}[t]{\subfigsizemed\textwidth}
			\centering
			\includegraphics[trim=0cm 0cm 0cm 0cm, clip, width=0.24\textheight]{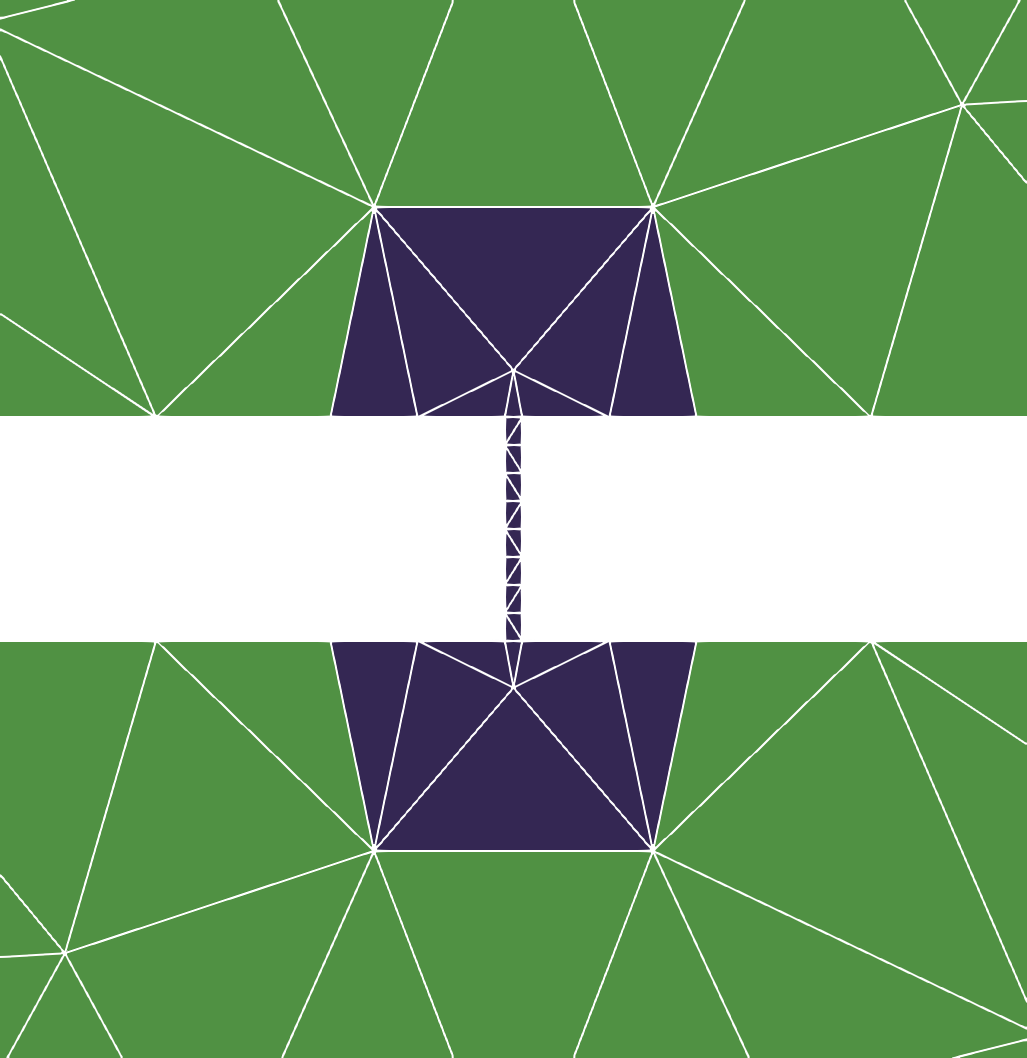}
			\caption{Zoom over the channel with $\delta=1/2^7$.}
			\label{fig:eff2_sto_gridb}
		\end{subfigure}
	\end{center}
	\caption{Narrow channel. Zoom of the FE mesh for channel width $\delta=1/2^2$ (left) or $1/2^7$ (right), with the subdomain $\Omega_{F,\delta}$ (in blue).}
	\label{fig:eff2_sto_grid}
\end{figure}

We illustrate the triangulation $\mathcal{M}$ in the neighborhood of the narrow channel in \cref{fig:eff2_sto_grid}, for two values of $\delta$. We observe that for large $\delta$ the typical element size is small enough to resolve the channel (\cref{fig:eff2_sto_grida}), while for small $\delta$ the elements in the channel are considerably smaller (\cref{fig:eff2_sto_gridb}). As the spectral radius of the discrete Laplacian behaves as $1/h^2$, where $h$ is the size of the smallest elements in the mesh, then $\rho$, the spectral radius of the Jacobian of $f$, increases as $\delta$ decreases. Therefore, the cost of SK-ROCK applied to \cref{eq:parode} increases as $\delta$ decreases.

Now, we want to decompose $f$ in two terms $\ff$ and $\fs$ such that as $\delta$ decreases then $\rhof$ increases but $\rhos$ remains constant, where $\rhof$ and $\rhos$ are the spectral radii of the Jacobians of $\ff$ and $\fs$, respectively.
We define a subdomain $\Omega_{F,\delta}$ consisting in the channel plus its neighboring elements having size smaller than the typical mesh size $H$, see \cref{fig:eff2_sto_grid}. Therefore, the size of the elements outside $\Omega_{F,\delta}$ is almost independent of $\delta$. In order to identify $\ff$ and $\fs$ as the discrete Laplacian inside and outside of $\Omega_{F,\delta}$, respectively, we define a diagonal matrix $D\in\Rb^{N\times N}$ by $D_{jj}=1$ if $\mbox{supp}(\varphi_j)\subset \Omega_{F,\delta}$ and $D_{jj}=0$ else. We let 
\begin{equation}
	\ff(X)=DAX, \qquad\qquad \fs(t,X)=(I-D)AX+M^{-1}\widehat b(t),
\end{equation}
with $I$ the identity matrix. Thus, as $\delta$ decreases, the size of the elements inside of $\Omega_{F,\delta}$ decrease and $\rhof$ increases, while $\rhos$ is independent of $\delta$.

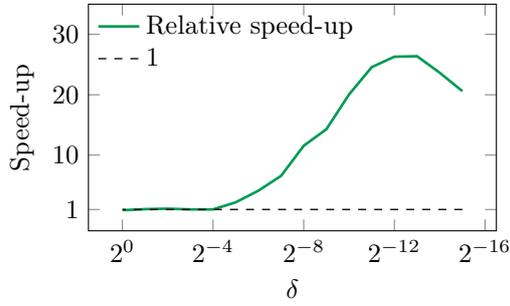
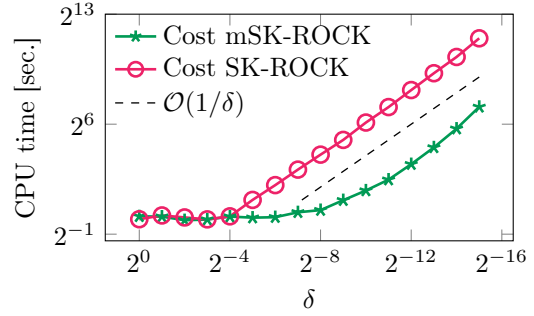
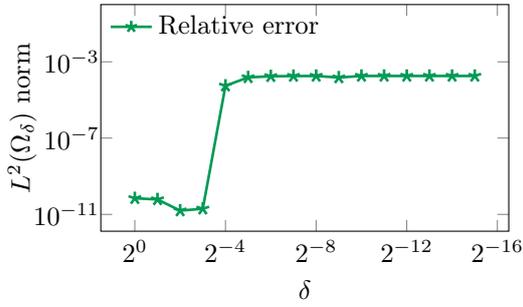
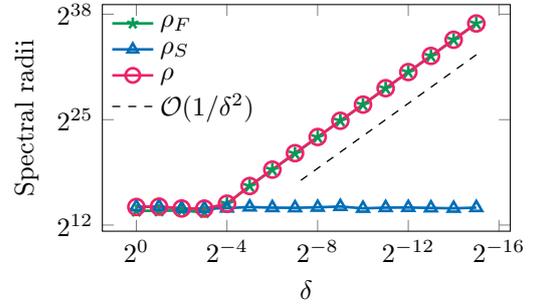
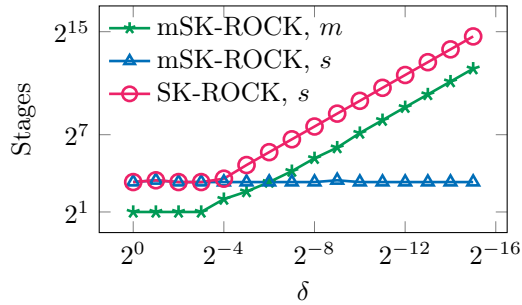
\begin{figure}
	\begin{center}
		\begin{subfigure}[t]{\subfigsize\textwidth}
			\centering
			\begin{tikzpicture}[scale=\plotimscale]
			\begin{semilogxaxis}[height=\aspectratio*\plotimsized\textwidth,width=\plotimsized\textwidth,legend columns=1,legend cell align={left},legend style={draw=\legendboxdraw,fill=\legendboxfill,at={(0,1)},anchor=north west},log basis x={2},x dir=reverse, ytick={1,10,20,30},ymax=35,
			xlabel={$\delta$}, ylabel={Speed-up},label style={font=\normalsize},tick label style={font=\normalsize},legend image post style={scale=\legendmarkscale},legend style={nodes={scale=\legendfontscale, transform shape}},]
			\addplot[color=colorone,solid,line width=\plotlinewidth pt,mark=none,mark size=\plotmarksizeu pt] table [x=delta,y=RelSpeedUp,col sep=comma] 
			{data/exp_eff2/results.csv};\addlegendentry{Relative speed-up}
			\addplot[black,dashed,line width=\plotdashedlinewidth pt,domain=3.0517578125e-05:1] (x,1);\addlegendentry{$1$}
			\end{semilogxaxis}
			\end{tikzpicture}
			\caption{Relative speed-up.}
			\label{fig:eff2_sto_speed}
		\end{subfigure}\hfill
		\begin{subfigure}[t]{\subfigsize\textwidth}
			\centering
			\begin{tikzpicture}[scale=\plotimscale]
			\begin{loglogaxis}[height=\aspectratio*\plotimsized\textwidth,width=\plotimsized\textwidth,legend columns=1,legend cell align={left},legend style={draw=\legendboxdraw,fill=\legendboxfill,at={(0,1)},anchor=north west},log basis x={2},log basis y={2},x dir=reverse,ymax=8192,ylabel style={yshift=-2pt},
			xlabel={$\delta$}, ylabel={CPU time [sec.]},label style={font=\normalsize},tick label style={font=\normalsize},legend image post style={scale=\legendmarkscale},legend style={nodes={scale=\legendfontscale, transform shape}},]
			\addplot[color=localcolor,solid,line width=\plotlinewidth pt,mark=\localmark,mark size=\plotmarksizeu pt] table [x=delta,y=time_mskrock,col sep=comma] 
			{data/exp_eff2/results.csv};\addlegendentry{Cost $\mSKROCK$}
			\addplot[color=classicalcolor,solid,line width=\plotlinewidth pt,mark=\classicalmark,mark size=\plotmarksized pt] table [x=delta,y=time_skrock,col sep=comma] 
			{data/exp_eff2/results.csv};\addlegendentry{Cost SK-ROCK}
			\addplot[black,dashed,line width=\plotdashedlinewidth pt,domain=3.0517578125e-05:0.0078] (x,1/x/64);\addlegendentry{$\bigo{1/\delta}$}
			\end{loglogaxis}
			\end{tikzpicture}
			\caption{Total CPU time w.r.t. $\delta$.}
			\label{fig:eff2_sto_comptime}
		\end{subfigure}\\ \vspace{0.5cm}
		\begin{subfigure}[t]{\subfigsize\textwidth}
			\centering
			\begin{tikzpicture}[scale=\plotimscale]
			\begin{loglogaxis}[height=\aspectratio*\plotimsized\textwidth,width=\plotimsized\textwidth,legend columns=1,legend cell align={left},legend style={draw=\legendboxdraw,fill=\legendboxfill,at={(0,1)},anchor=north west},log basis x={2},log basis y={10},x dir=reverse, ytick={1e-11,1e-7,1e-3},ymax=1,ylabel style={yshift=-8pt},yticklabel style={xshift=3pt},
			xlabel={$\delta$}, ylabel={$L^2(\Omega_\delta)$ norm},label style={font=\normalsize},tick label style={font=\normalsize},legend image post style={scale=\legendmarkscale},legend style={nodes={scale=\legendfontscale, transform shape}},]
			\addplot[color=colorone,solid,line width=\plotlinewidth pt,mark=\markone,mark size=\plotmarksizeu pt] table [x=delta,y=relerrl2,col sep=comma] 
			{data/exp_eff2/results.csv};\addlegendentry{Relative error}%{$\Vert y_N^{\RKCop}-y_N^{\mRKCop}\Vert_{H^1(\Omega_\delta)}$}
			%\addplot[black,dashed,domain=6.103515625e-05:1] (x,0.001);\addlegendentry{$10^{-3}$}
			\end{loglogaxis}
			\end{tikzpicture}
			\caption{SK-ROCK and $\mSKROCK$ solutions' relative error $\Vert u^{\mSKROCKop}-u^{\SKROCKop}\Vert/\Vert u^{\SKROCKop}\Vert$ in $L^2(\Omega_\delta)$ norm.}
			\label{fig:eff2_sto_relerr}
		\end{subfigure}\hfill
		\begin{subfigure}[t]{\subfigsize\textwidth}
			\centering
			\begin{tikzpicture}[scale=\plotimscale]
			\begin{loglogaxis}[height=\aspectratio*\plotimsized\textwidth,width=\plotimsized\textwidth,legend columns=1,legend style={draw=\legendboxdraw,fill=\legendboxfill,at={(0,1)},anchor=north west},log basis x={2},log basis y={2},x dir=reverse,legend cell align={left},
			xlabel={$\delta$}, ylabel={Spectral radii},label style={font=\normalsize},tick label style={font=\normalsize},legend image post style={scale=\legendmarkscale},legend style={nodes={scale=\legendfontscale, transform shape}},]
			\addplot[color=localcolorF,solid,line width=\plotlinewidth pt,mark=\localmarkF,mark size=\plotmarksizeu pt] table [x=delta,y=rhoF_mskrock,col sep=comma] 
			{data/exp_eff2/results.csv};\addlegendentry{$\rhof$}
			\addplot[color=localcolorS,solid,line width=\plotlinewidth pt,mark=\localmarkS,mark size=\plotmarksizeu pt] table [x=delta,y=rhoS_mskrock,col sep=comma] 
			{data/exp_eff2/results.csv};\addlegendentry{$\rhos$}
			\addplot[color=classicalcolor,solid,line width=\plotlinewidth pt,mark=\classicalmark,mark size=\plotmarksized pt] table [x=delta,y=rho_skrock,col sep=comma] 
			{data/exp_eff2/results.csv};\addlegendentry{$\rho$}
			\addplot[black,dashed,line width=\plotdashedlinewidth pt,domain=3.0517578125e-05:0.0078] (x,8/x/x);\addlegendentry{$\bigo{1/\delta^2}$}
			\end{loglogaxis}
			\end{tikzpicture}
			\caption{Spectral radii w.r.t. $\delta$.}
			\label{fig:eff2_sto_rho}
		\end{subfigure}\\ \vspace{0.5cm}
		\begin{subfigure}[t]{\subfigsize\textwidth}
			\centering
			\begin{tikzpicture}[scale=\plotimscale]
			\begin{loglogaxis}[height=\aspectratio*\plotimsized\textwidth,width=\plotimsized\textwidth,legend columns=1,legend style={draw=\legendboxdraw,fill=\legendboxfill,at={(0,1)},anchor=north west},log basis x={2},log basis y={2},x dir=reverse,legend cell align={left},legend cell align={left},ytick={2,128,32768},ymax=131072,
			xlabel={$\delta$}, ylabel={Stages},label style={font=\normalsize},tick label style={font=\normalsize},legend image post style={scale=\legendmarkscale},legend style={nodes={scale=\legendfontscale, transform shape}},]
			\addplot[color=localcolorF,solid,line width=\plotlinewidth pt,mark=\localmarkF,mark size=\plotmarksizeu pt] table [x=delta,y=m_mskrock,col sep=comma] 
			{data/exp_eff2/results.csv};\addlegendentry{$\mSKROCK$, $m$}
			\addplot[color=localcolorS,solid,line width=\plotlinewidth pt,mark=\localmarkS,mark size=\plotmarksizeu pt] table [x=delta,y=s_mskrock,col sep=comma] 
			{data/exp_eff2/results.csv};\addlegendentry{$\mSKROCK$, $s$}
			\addplot[color=classicalcolor,solid,line width=\plotlinewidth pt,mark=\classicalmark,mark size=\plotmarksized pt] table [x=delta,y=s_skrock,col sep=comma] 
			{data/exp_eff2/results.csv};\addlegendentry{SK-ROCK, $s$}
			\end{loglogaxis}
			\end{tikzpicture}
			\caption{Number of stages.}
			\label{fig:eff2_sto_stages}
		\end{subfigure} 
		%\begin{subfigure}[t]{\subfigsize\textwidth}
		%\centering
		%\begin{tikzpicture}[scale=\plotimscale]
		%\begin{semilogxaxis}[height=\aspectratio*\plotimsized\textwidth,width=\plotimsized\textwidth,legend columns=1,legend cell align={left},legend style={draw=\legendboxdraw,fill=\legendboxfill,at={(0,1)},anchor=north west},log basis x={2},x dir=reverse, ytick={1,20,40,60},ymax=75,
		%xlabel={$\delta$}, ylabel={Speed-up},label style={font=\normalsize},tick label style={font=\normalsize},legend image post style={scale=\legendmarkscale},legend style={nodes={scale=\legendfontscale, transform shape}},]
		%\addplot[color=colorone,solid,line width=\plotlinewidth pt,mark=none,mark size=\plotmarksizeu pt] table [x=delta,y=eS,col sep=comma] 
		%{mSKROCK/data/exp_eff2/results_nonoise.csv};\addlegendentry{Speed-up: no noise case}
		%\addplot[black,dashed,domain=3.0517578125e-05:1] (x,1);\addlegendentry{$1$}
		%\end{semilogxaxis}
		%\end{tikzpicture}
		%\caption{Relative speed-up in the deterministic case.}
		%\label{fig:eff2_sto_speed_nonoise}
		%\end{subfigure}
	\end{center}
	\caption{Narrow channel. Efficiency, evolution of spectral radii and stages.}
	\label{fig:eff2_sto_results}
\end{figure}

We will solve \cref{eq:parode} for varying channel width $\delta$ and investigate the efficiency of the $\mSKROCK$ and SK-ROCK method. Hence, for each $\delta=1/2^k$ with $k=0,\ldots,15$ we solve once
\begin{equation}
	\dif X(t)=\fs(t,X(t))\dif t+\ff(X(t))\dif t+g(X(t))\dif B(t)\quad t\in (0,T], \qquad\qquad X(0)=0
\end{equation}
with the $\mSKROCK$ and SK-ROCK methods, on the same sample path $B(t)$ with $T=0.1$ and the same step size $\tau=0.01$. The relative speed-up $S$ given by the $\mSKROCK$ scheme over the SK-ROCK method, in terms of CPU time, in function of $\delta$ is displayed in \cref{fig:eff2_sto_speed}. For large $\delta$ both methods have the same performance ($S\approx 1$), as $\delta$ decreases the $\mSKROCK$ becomes more efficient than SK-ROCK and it is at least 25 times faster for some values of $\delta$. 

The relative speed-up has been computed dividing the computational costs (CPU time) of the SK-ROCK and $\mSKROCK$ method, that are plotted in \cref{fig:eff2_sto_comptime}. This choice is justified by the fact that the relative error between the two solutions, measured in the $L^2(\Omega_\delta)$ norm at time $T$, is less than 1\textperthousand, see \cref{fig:eff2_sto_relerr}. 
Note that in \cref{fig:eff2_sto_relerr} a jump appears exactly when $m$ passes from $m=2$ to $m\geq 4$ (see \cref{fig:eff2_sto_stages}) and thus $r$ passes from $r=1$ to $r\geq 2$. This is due to the fact that for smaller $m$ the $\mSKROCK$ and SK-ROCK schemes are closer.

The spectral radii $\rho,\rhof,\rhos$ of the Jacobians of $f,\ff,\fs$ are shown in \cref{fig:eff2_sto_rho}, for large $\delta$ the typical element size is sufficiently small to resolve the channel (\cref{fig:eff2_sto_grida}) and thus $\rho\approx\rhof\approx\rhos$, implying that the costs of $\mSKROCK$ and SK-ROCK are similar. As $\delta$ decreases then $\rho,\rhof$ increase. Since $\rhos$ is almost constant the number of $\fs$ evaluations in the $\mSKROCK$ method remains constant and only the number of $\ff$ evaluations increase, therefore the cost of $\mSKROCK$ increases less rapidly than the one of SK-ROCK. Finally, in \cref{fig:eff2_sto_stages} we show the number of stages taken by the methods, which reflects the behavior of the spectral radii. 

In \cref{fig:eff2_sto_speed} we see a decrease in speed-up for $\delta$ extremely small, this is due to the fact that the cost of evaluating $\ff$, with respect to $\fs$ and $g$, becomes important; a high number of tiny elements is indeed needed to resolve the channel, hence the total cost $c_F+c_S+c_g=1$ is dominated by $c_F$ close to $1$ (see \cref{sec:mskrockalg}) and this is not the optimal speed-up for the $\mSKROCK$ method. Nevertheless, the $\mSKROCK$ scheme still remains about 20 times faster than the SK-ROCK method. We note that in practical applications the mesh outside the channel would also be refined and a value of $c_F$ closer to the optimal speed-up could be reached.

\section{Conclusion}
We have introduced a modified equation 
%\cref{eq:modsde} 
$\dif \Xe(t) = \fe(X(t))\dif t +\ge(X(t))\dif W(t)$
%defined by an average force $\fe$ and a damped diffusion $\ge$ ,
for stiff stochastic differential equations 
%\cref{eq:msde}
$\dif X(t) = \ff(X(t))\dif t +\fs(X(t))\dif t+ g(X(t))\dif W(t)$
with different time-scales but without any clear-cut scale separation, where the drift is composed by a stiff but cheap term $\ff$ and a mildly stiff but expensive term $\fs$. 
%As there is no scale separation, $\ff$ might contain slow terms too.
The averaged force $\fe$ is such that the stiffness of the modified equation depends solely on the slow term $\fs$, while the damped diffusion $\ge$ is such that the mean-square stability properties of the original problem are preserved. Therefore, integration of the modified equation by explicit schemes is cheaper than the original problem, as the stability conditions are not affected by a few severely stiff degrees of freedom in $\ff$. Evaluation of both $\fe,\ge$ requires the solution to fast but cheap deterministic auxiliary problems \cref{eq:defu,eq:defv}, which can be approximated by explicit schemes.

Starting from the modified equation we devised an interpolation-free stabilized explicit multirate scheme, given by \cref{eq:mskrock,eq:defbge,eq:defve,eq:defbve,eq:defbfe,eq:defbue}. The method consists in integrating the modified equation with a stabilized explicit scheme for SDEs (SK-ROCK) and evaluating $\fe,\ge$ by solving the auxiliary problems with a stabilized explicit scheme for ODEs (RKC). 
The scheme, called $\mSKROCK$, is fully explicit, has strong order $1/2$, weak order $1$ and \modr{it is proven to be stable on a model problem} --- see \cref{thm:stab_mskrock,thm:conv}. The number of expensive function evaluations of $\fs$ needed by $\mSKROCK$ depends only on $\fs$ itself; therefore, the efficiency of the scheme is hardly affected by the severely stiff term $\ff$.

Furthermore, an important property of the scheme is that it is not based on any scale separation assumption. Therefore, it can be employed for systems stemming from the spatial discretization of stochastic parabolic partial differential equations on locally refined grids, where $\ff,\fs$ represent the Laplacian in refined and coarse regions, respectively (see \cref{sec:channel}).
Finally, the method is straightforward to implement and numerical experiments demonstrate that the computational cost is significantly reduced without sacrificing any accuracy, compared to the optimal stabilized method for stiff SDEs, namely the SK-ROCK method.

\section*{Acknowledgments} The authors are partially supported by the Swiss National Science Foundation, under grant No. $200020\_172710$.

%\bibliographystyle{habbrv}
%\bibliographystyle{abbrv}
%\bibliography{../../../../../../LaTeX/library}

\end{document}